\documentclass{article}
\usepackage{amsmath}
\usepackage{amssymb}
\usepackage{amsthm}
\usepackage{amscd}
\newtheorem{thm}{Theorem}[section]
\newtheorem{prop}[thm]{Proposition}
\newtheorem{lem}[thm]{Lemma}
\newtheorem{cor}[thm]{Corollary}
\newtheorem{defn}[thm]{Definition}
\theoremstyle{remark}
\newtheorem{rem}{Remark}

\newcommand{\al}{\alpha}
\newcommand{\be}{\beta}
\newcommand{\ga}{\gamma}
\newcommand{\de}{\delta}

\newcommand{\mbb}{\mathbb}
\newcommand{\fra}{\mathfrak}
\newcommand{\ov}{\overline}
\newcommand{\vphi}{\varphi}

\begin{document}

\title{Groups Elementarily Equivalent to a Free 2-nilpotent Group of Finite Rank}
\author{Alexei G. Myasnikov \and Mahmood Sohrabi}
\maketitle

\section*{Abstract}In this paper we find a characterization for groups elementarily equivalent to a free nilpotent group $G$ of class 2 and arbitrary finite rank.

\tableofcontents

\section{Introduction}

In this paper we give an algebraic description of groups elementarily equivalent to a given free nilpotent class 2 group of an arbitrary finite rank.

\subsection{Elementary classification problem in groups}

Importance of the elementary (logical) classification of algebraic structures was emphasized in the works of A.Tarski and A.Malcev. In general, the problem of elementary classification requires to characterize, in somewhat algebraic terms, all algebraic structures (perhaps, from a given class) elementarily equivalent to a given one. Recall, that two algebraic structures $\mathcal{A}$ and $\mathcal{B}$ in a language $L$ are elementarily equivalent ($\mathcal{A} \equiv \mathcal{B}$) if they have the same first-order theories in $L$ (undistinguishable in the first order logic in $L$).

Tarski's results \cite{T51} on elementary theories of $\mathbb{R}$ and $\mathbb{C}$ (algebraically closed and  real closed fields), as well as, the subsequent results of J. Ax and S. Kochen \cite{AK1, AK2,AK3} and Yu. Ershov \cite{ershov1, ershov2, ershov3} on elementary theories  of $\mathbb{Q}_p$ ($p$-adically closed fields) became an algebraic classic and can be found in many text-books on model theory.

One of the initial influential results on elementary classification of groups is due to W.Szmielew - she classified elementary theories of abelian groups in terms of "Szmielew's" invariants \cite{Szmielew} (see also \cite{eklof,Monk, Baur}).  For non-abelian groups, the main inspiration, perhaps,  was  the famous Tarski's problem whether free non-abelian groups of finite rank are elementarily equivalent or not. This problem has been open for many years and only recently was solved in affirmative in \cite{KM3,Sela6}. In contrast, free solvable (or nilpotent) groups of finite rank are elementarily equivalent if and only if they are isomorphic (A. Malcev \cite{Malcev1}).  Indeed, in these cases the abelianization $G/[G,G]$ of the group $G$ (hence the rank of $G$) is definable  (interpretable)   in $G$  by first-order formulas, hence the result.

In  \cite{malcev3} A. Malcev described elementary equivalence among classical linear groups. Namely, he showed that if $\mathcal{G} \in \{GL,PGL,SL,PSL\}$, $n, m \geq 3$, $K$ and $F$ are fields of characteristic zero, then $\mathcal{G}(F)_m \equiv \mathcal{G}(K)_n$ if and only if $m = n$ and $F \equiv K$. It turned out later that  this type of  results can be obtained via   ultrapowers by means of  the theory of abstract isomorphisms of such groups.   In this approach one argues that if the groups $\mathcal{G}(F)_m$ and $\mathcal{G}(K)_n$  are elementarily equivalent then their ultrapowers over a  non-principal ultrafilter $\omega$ are isomorphic. Since these  ultrapowers are again  groups of the type $\mathcal{G}(F^*)_m$ and $\mathcal{G}(K^*)_n$ (where $F^*$ and $K^*$ are the corresponding ultrapowers of the fields) the result follows from the description of abstract isomorphisms of such groups (which are semi-algebraic, so they preserve the algebraic scheme and the field). It follows that the results, similar to the ones mentioned above,   hold for many algebraic and linear groups. We refer to  \cite{RR} and \cite{BM04,BM05} for details. On the other hand, many "geometric" properties of algebraic groups are just first-order definable invariants of these groups, viewed as abstract groups (no  geometry, only multiplication). For example, the geometry of a simple algebraic group is entirely determined by its group multiplication (see \cite{Z84,Poi88,Poibook}), which readily implies the celebrated Borel-Tits theorem on abstract isomorphisms of simple algebraic groups.

Other large classes  of groups where the elementary classification  problem is relatively well-understood are the classes of finitely generated nilpotent groups and algebraic nilpotent groups. This is the main subject  of the paper and we discuss it below. We finish a short survey on elementary equivalence  of groups in a   result  related to the solutions of the Tarski's problem. It is known now  that a finitely generated group $G$ is elementarily equivalent to a free non-abelian group if and only if $G$ is a regular NTQ-group (\cite{KM3}), which are  the same as $\omega$-residually free towers (\cite{Sela6}). These  results generated a new waive of research on elementary classification in classes of groups which are very different from solvable or algebraic groups, namely, in hyperbolic groups, relatively hyperbolic groups, free products of groups, right angled Artin groups, etc. On the other hand, research on large scale geometry of soluble groups, in particular, their quasi-isometric invariants, comes suspiciously close to the study of their first-order invariants. That gives a new impetus to learn more on elementarily equivalence  in the class of solvable groups.

\subsection{On elementary classification of nilpotent groups}

There are several principal results known  on elementary theories of nilpotent groups. In his pioneering paper \cite{malcev2} A. Malcev  showed that a ring $R$ with unit can be defined by first-order formulas in the group $UT_3(R)$ of unitriangular matrices over $R$ (viewed as an abstract group). In particular, the ring of integers $\mathbb{Z}$ is definable in the group $UT_3(\mathbb{Z})$, which is a free 2-nilpotent group of rank 2. In \cite{ershov4} Yu. Ershov  proved that the group $UT_3(\mathbb{Z})$ (hence the ring $\mathbb{Z}$) is definable in any finitely generated nilpotent group $G$, which is not virtually abelian. It follows immediately that the elementary theory of $G$ is undecidable.
On the elementary classification side the main research was on M.Kargapolov's conjecture: two finitely generated nilpotent groups are elementarily equivalent if and only if they are isomorphic. In \cite{Z71}  B. Zilber gave a counterexample to the Kargapolov's conjecture. In the break-through papers \cite{MR1,MR2,MR3} A. Myasnikov and V. Remeslennikov proved that the Kargapolov's conjecture holds "essentially" true in the class of  nilpotent $\mathbb{Q}$-groups (i.e., divisible torsion-free nilpotent groups) finitely generated as $\mathbb{Q}$-groups. Indeed, it turned out that two such groups $G$ and $H$ are elementarily equivalent if their {\it cores} $\bar{G}$ and $\bar{H}$ are isomorphic and $G$ and $H$ either simultaneously  coincide with their cores or they do not. Here the core  of $G$ is uniquely defined as a subgroup $\bar{G} \leq G$ such that $Z(\bar{G}) \leq [\bar{G},\bar{G}]$ and $G= \bar{G} \times G_0$,  for some  abelian $\mathbb{Q}$-group  $G_0$. Developing this approach further A.Myasnikov described in \cite{alexei87,alexei90a} all groups elementarily equivalent to a given finitely generated nilpotent $K$-group $G$ over an arbitrary field of characteristic zero. Here by a $K$-group we understand P. Hall nilpotent $K$-powered groups, which are the same as $K$-points of nilpotent algebraic groups, or unipotent $K$-groups. Again, the crucial point here is that the geometric structure of the group $G$ (including the fields of definitions of the components of $G$ and their related structural constants) are first-order definable in $G$, viewed as an abstract group. Furthermore, these ideas shed some light on the Kargapolov's conjecture - it followed that two finitely generated elementarily equivalent nilpotent groups $G$ and $H$ are isomorphic, provided one of them is a {\it core} group. In this case  $G$ is a {\em core} group if $Z(G) \leq I([G,G])$, where $I([G,G])$ is the isolator of the commutant $[G,G]$. Finally, F.Oger showed in \cite{oger} that two finitely generated nilpotent groups $G$ and $H$ are elementarily equivalent if and only if they are essentially isomorphic, i.e., $G \times \mathbb{Z} \simeq H \times \mathbb{Z}$.
However, the full  classification problems for finitely generated nilpotent groups is currently wide open. In a series of papers \cite{beleg92,beleg94,beleg99} O. Belegradek completely characterized groups which are elementarily equivalent to a nilpotent group $UT_n(\mathbb{Z})$ for $n \geq 3$. It  is easy to see that (via ultrapowers) that if $\mathbb{Z} \equiv R$ for some ring $R$ then $UT_n(\mathbb{Z}) = UT_n(R)$. However, it has been shown in \cite{beleg94,beleg99} that there are groups   elementarily equivalent to $UT_n(\mathbb{Z})$ which  are not isomorphic to any group of the type $UT_n(R)$ as above (quasi-unitriangular groups).

\subsection{Results and the structure of the paper}

In this paper we generalize O.  Belegradek's results  from ~\cite{beleg92} on characterizing groups elementary equivalent to the group $UT_3(\mathbb{Z})$, which is, as mentioned above, a free nilpotent group of class 2 and rank 2. We start with describing a certain class $N_{2,n}$ of groups such that $G=N_{2,n}(\mathbb{Z})$ is a free 2-nilpotent group of rank $n$ for a natural number $n\geq 2$. The class $N_{2,2}$ coincides with the class of upper unitriangular groups over a commutative associative unitary rings $R$. We introduce various definitions for a ``basis'' of an $N_{2,n}$ group over $R$  and show that these bases are first order definable.  As should be expected from the work of O. Belegradek on $UT_3$ groups the classes $N_{2,n}$ are not elementarily closed. To close the class   we introduce a new  type of groups $QN_{2,n}$, and prove that they give   the elementary closure of the class of  $N_{2,n}$ groups over commutative unitary rings. One of the most crucial steps in this work is ``recovering'' the ring of integers $\mathbb{Z}$ from $G$ and showing that it is absolutely interpretable in $G$. To do so we use the method of bilinear mappings, due to A. Myasnikov \cite{alexei86,alexei91}, in contrast with O. Belegradek who uses Mal'cev original idea from~\cite{malcev2}.

Now we describe contents of the sections. In Section~\ref{prel} we present the basic facts and definitions from
the theory of nilpotent groups, the theory of group extensions, model theory and model theory of bilinear mappings. The classes of free 2-nilpotent groups $N_{2,n}$ of rank $n$ over commutative associative rings $R$ with unit and some related concepts are discussed  in Section~\ref{n2ns}. In Section~\ref{quasiqnn2}  we introduce  a new class of nilpotent groups $QN_{2,n}$, termed quasi  $N_{2,n}$ groups.
 We give a characterization theorem for the
$QN_{2,n}$ groups over $R$ Section~\ref{charthmsec}. This section contains the main result, Theorem~\ref{main}, of this paper in which we prove that a group elementarily equivalent to a
free 2-nilpotent group of arbitrary rank $n$ is a group of the type $QN_{2,n}$ over a ring $R$ with $R \equiv \mathbb{Z}$. In Section~\ref{cse} we prove that the necessary condition above is also a sufficient condition, thus providing a complete characterization of groups elementarily equivalent to a given free 2-nilpotent  groups of finite rank.
In Section \ref{qnnotn} we show an example of a group  $QN_{2,n}$, which is elementarily equivalent to $N_{2,n}(\mathbb{Z})$ but is not an $N_{2,n}$ group over any commutative associative unitary ring.

\section{Preliminaries and notation}\label{prel}
\subsection{Nilpotent groups}

\label{upperlower}

Let $G$ be a group with a series of subgroups:
$$G=G_1\geq G_2 \geq \ldots G_n \geq G_{n+1}=1,$$
where each $G_{i+1}$ is a normal subgroup of $G$ and each factor $G_i/G_{i+1}$ is an abelian group. Let $G$
act on each factor $G_i/G_{i+1}$ by conjugation, i.e.
$$g.xG_{i+1}=_{df}g^{-1}xgG_{i+1}.$$
If the above action of $G$ on all the factors is trivial then the above series is called a \textit{central
series} and any group $G$ with such a series is called a \textit{nilpotent} group.

 For elements $x$ and $y$ of a group $G$ let $[x,y]=x^{-1}y^{-1}xy$. $[x,y]$ is
called the \textit{commutator} of the elements $x$ and $y$. The subgroup $[G,G]$ is the subgroup of $G$
generated by all $[x,y]$, $x,y\in G$. In general for $H$ and $K$ subgroups of $G$, $[H,K]$ is the subgroup of
$G$ generated by commutators $[x,y]$, $x\in H$ and $y\in K$. Let us define a series $\Gamma_1(G), \Gamma_2(G),
\ldots$ of subgroups of $G$ by setting
$$G=\Gamma_1(G), \quad \Gamma_{n+1}(G)=[\Gamma_n(G),G]\quad \textrm{for all $n>1$}.$$
It can be easily checked that the above series is a central series. If $c$ is the least number that
$\Gamma_{c+1}(G)=0$ then $G$ is said to be a nilpotent group of \textit{class $c$} or simply a
\textit{$c$-nilpotent} group. We call the series above the \textit{lower central series} of the group $G$.

Let $Z(G)$ denote the center of a group $G$. We define a series of subgroups $Z_i(G)$ of $G$ by setting
$$Z_1(G)=Z(G), \quad Z_{i+1}(G)=\{x \in G:xZ_i\in Z(G/Z_i(G))\},\quad i> 1.$$ This series is also a central series and called the \textit{upper central
series} of the group $G$. If $Z_{n+1}(G)=G$ for some finite number $n$ and $c$ is the least such number then $G$
is provably a $c$-nilpotent group.

Let $F(n)$ be the free group on $n$ generators. Let $G$ be a group isomorphic to the factor group
$F(n)/\Gamma_{c+1}(F(n))$. Then $G$ is called a \textit{free nilpotent group of class $c$ and rank n}, or equivalently a \textit{free c-nilpotent group of rank n}.

For definitions and details regarding Mal'cev (canonical) basis of a finitely generated torsion free nilpotent group, groups admitting exponents in a binomial domain $R$ or $R$-group for short and $R$-completions of finitely generated torsion free nilpotent groups we refer to~\cite{hall}.

\subsection{\protect Central and abelian
extensions}\label{extensions}
 Let $A$ and $B$ be abelian groups. Consider the short exact sequence:
$$0\rightarrow A \xrightarrow{\mu} E \xrightarrow{\nu} B \rightarrow 0.$$
 The group $E$ is called an \textit{abelian extension} of $A$ by $B$ if $E$ is an abelian group. $E$ is said to be a
\textit{central extension} of $A$ by $B$ if $\mu(A)$ sits inside the center of $E$, i.e. the action defined
above is trivial. Obviously every abelian extension is central. It can be easily seen that every central
extension of two abelian groups is a 2-nilpotent group. An extension:
$$0\rightarrow A \xrightarrow{\mu'} E' \xrightarrow{\nu'} B \rightarrow 0$$
is \textit{equivalent} to the extension above if there is an isomorphism $\eta:E \rightarrow E'$ such that
$\nu'\circ \eta= \nu$ and $\eta \circ \mu= \mu'$. The relation ``equivalence'' defines an equivalence relation
on the set of all central extensions of the abelian groups $A$ and $B$.

We now review the relation between equivalence classes of central extensions of an abelian group $A$ by an
abelian group $B$ and the group called the \textit{second cohomology} group, $H^2(B,A)$.

A function $f:B\times B\rightarrow A$  is called a \textit{2-cocycle} if
\begin{itemize}
\item $f(0,x)=f(x,0)=0$ for every $x$ in $B$,
\item $f(x+y,z) + f(x,y)=f(x,y+z) +f(y,z)$, for $x$, $y$ and $z$ in $B$.
\end{itemize}
By $B^2(B,A)$ we denote the set of all 2-cocycles $f:B\times B\rightarrow A$. A We can make the set $B^2(B,A)$ into
 an abelian group by letting addition of the corresponding functions be the point-wise addition. An element $f\in B^2(B,A)$ is called a \textit{2-coboundary} if there exists a function $\psi:B\rightarrow A$ such that
 $$f(x,y)=\psi(x+y)-\psi(x)-\psi(y).$$
 We note the set of all 2-coboundaries by $Z^2(B,A)$. The group operation defined above on $B^2(B,A)$ makes $Z^2(B,A)$ into a subgroup. An element $f\in B^2(B,A)$ is called a \textit{symmetric 2-cocycle} if
 $$f(x,y)=f(y,x), \quad \forall x,y \in B.$$
By $S^2(B,A)$ we mean the subgroup of $B^2(B,A)$ consisting of all symmetric 2-cocycles. Set
$$H^2(B,A)=_{df}B^2(B,A)/Z^2(B,A),$$
and
$$Ext(B,A)=_{df}S^2(B,A)/(Z^2(B,A)\cap S^2(B,A).$$
There is a 1-1 correspondence between elements of $H^2(B,A)$ and equivalence classes of central extensions of $A$ by $B$. The same correspondence also exists between elements of $Ext(B,A)$ and equivalences classes of abelian extensions of $A$ by $B$. We briefly explain one direction the correspondence. For more details we refer to~\cite{robin}. let $f:B\times B\rightarrow A$ be a 2-cocycle.
Define a group $E(f)$ by $E(f)=B\times A$ as sets with the multiplication
$$(b_1, a_1)(b_2, a_2)=(b_1+b_2, a_1+a_2+f(b_1,b_2)) \quad a_1,a_2\in A, \quad b_1,b_2 \in B.$$
The above operation is a group operation and the resulting extension is central. If $f$ is symmetric then $E(f)$ is abelian. Moreover it can be verified
that if $f,f':B^(B,A)$ and $f-f'\in Z^2(B,A)$ then the extensions
$E(f)$ and $E(f')$ are equivalent.

 \subsection{\protect Structures, signatures and interpretations}\label{ssi}
\label{interpret1}

A group $G$ is considered to be
the structure $\langle|G|,.,^{-1},1\rangle$ where $.$, $^{-1}$ and $1$, name multiplication, inverse operation
and the trivial element of the group respectively. We consider this signature as the \textit{signature of
groups}. We use $[x,y]$ as an abbreviation for $x^{-1}.y^{-1}.x.y$.

By an \textit{algebraic structure} we mean a structure including functions only, constants aside. Strangely
enough here we assume that algebraic structures consist only of predicates in addition to
constant symbols. But in a sense what we mean is clear. Algebraic operations are considered as relations rather
than functions.

Let $\mathfrak{U}$ be a structure and $\phi(x_1,\ldots, x_n)$ be a first order formula of the signature of
$\mathfrak{U}$ with $x_1$,\ldots ,$x_n$ free variables. Let $(a_1,\ldots, a_n)\in |\mathfrak{U}|^n$. We denote
such a tuple by $\bar{a}$. The notation $\mathfrak{U} \models \phi(\bar{a})$ is intended to mean that the tuple
$\bar{a}$ satisfies $\phi(\bar{x})$ when $\bar{x}$ is an abbreviation for the tuple $(x_1,\ldots, x_n)$ of
variables. For definitions of a formula of a signature, free variables and satisfaction the reader
should refer to \cite{hodges}.

Given a structure $\mathfrak{U}$ and a first order formula $\phi(x_1,\ldots, x_n)$ of the signature of
$\mathfrak{U}$, $\phi(\mathfrak{U}^n)$ refers to $\{\bar{a}\in |\mathfrak{U}|^n: \mathfrak{U}\models
\phi(\bar{a})\}$. Such a relation or set is called \textit{first order definable without parameters}. If
$\psi(x_1,\ldots,x_n,y_1,\ldots,y_m)$ is a first order formula of the signature of $\mathfrak{U}$ and $\bar{b}$
an $m$-tuple of elements of $\mathfrak{U}$ then $\psi(\mathfrak{U}^n,\bar{b})$ means $\{\bar{a}\in
|\mathfrak{U}|^n :\mathfrak{U}\models \psi(\bar{a},\bar{b})\}$. A set or relation like this is said to be
\textit{first order definable with parameters}.

Let $\mathfrak{U}$ be a structure of signature $\Sigma$. The \textit{theory} $Th(\mathfrak{U})$ of the structure
$\mathfrak{U}$ is the set: $$\{\phi:\mathfrak{U}\models \phi, \phi \textrm{ a first order sentence of signature
}\Sigma\}.$$

Finally two structures $\mathfrak{U}$ and $\mathfrak{B}$ of the signature $\Sigma$ are \textit{elementarily
equivalent} if $Th(\mathfrak{U})=Th(\mathfrak{B})$.

Let $\mathfrak{B}$ and $\mathfrak{U}$ be algebraic structures of signatures $\Delta$ and
$\Sigma$ respectively not having function symbols. In the following $\bar{x}$ refers to an $n$-tuple of variables $(x_1, \ldots , x_n)$ and for each $i$, $\ov{y^i}$ refers to an $m$-tuple of variables $(y^i_1, \ldots , y^i_m)$. The structure $\mathfrak{U}$ is said to be
\textit{interpretable} in $\mathfrak{B}$ with parameters $\bar{b}\in |\mathfrak{B}|^n$ or \textit{relatively
interpretable} in $\mathfrak{B}$ if there is a set of first order formulas
$$\Psi=\{A(\bar{x},\bar{y}), E(\bar{x},\ov{y^1},\ov{y^2}),\psi_{\sigma}(\bar{x}, \ov{y^1}, \ldots ,
\ov{y^{t_{\sigma}}}): \sigma \textrm{  a predicate of signature  } \Sigma \}$$ of the signature $\Delta$ such
that
 \begin{enumerate}
 \item $A(\bar{b})=\{\bar{a}\in|\mathfrak{B}|^m:\mathfrak{B}\models A(\bar{b},\bar{a})\}$ is not empty,
 \item $E(\bar{x},\ov{y^1},\ov{y^2})$ defines an equivalence relation $\epsilon_{\bar{b}}$ on $A(\bar{b})$,
 \item if the equivalent class of a tuple of elements $\bar{a}$ from $A(\bar{b})$ modulo the
 equivalence relation $\epsilon_{\bar{b}}$ is denoted by $[\bar{a}]$, for every predicate
 $\sigma$ of signature $\Sigma$, predicate $P_{\sigma}$ is defined on
 $A(\bar{b})/\epsilon_{\bar{b}}$ by
 $$P_{\sigma}([\bar{b}],[\ov{a^1}], \ldots [\ov{a^{t_\sigma}}])\Leftrightarrow_{df}\mathfrak{B}\models \psi_{\sigma}(\bar{b}, \ov{a^1},
 \ldots \ov{a^{t_\sigma}}),$$
 \item the structures $\mathfrak{U}$ and
 $\Psi(\mathfrak{B},\bar{b})=\langle A(\bar{b})/\epsilon_{\bar{b}},P_{\sigma}:\sigma\in \Sigma \rangle$ are isomorphic.
 \end{enumerate}
 Let $\phi(x_1,\ldots, x_n)$ be a first order formula of the signature $\Delta$ and $\bar{b}\in \phi(\mathfrak{B}^n)$ be as above. If $\mathfrak{U}$ is interpretable in $\mathfrak{B}$ with the
 parameters $\bar{b}$ and $\mathfrak{B}\models \phi(\bar{b})$ then  $\mathfrak{U}$
 is said to be \textit{regularly interpretable} in $\mathfrak{B}$ with the help of formula $\phi$. If the tuple $\bar{b}$ is empty, $\mathfrak{U}$ is said
 be \textit{absolutely interpretable} in $\mathfrak{B}$.

\subsection{\protect Bilinear mappings as model theoretic objects}\label{modelbilin}
All the results in this subsection are due to A. G. Myasnikov, \cite{alexei86} and \cite{alexei91}.

Let $M$ and $N$ be exact $R$-modules for some commutative ring $R$. An $R$-module $M$ is \textit{exact} if
$rm=0$ for $r\in R$ and $0\neq m\in M$ imply $r=0$. Let's recall that an $R$-bilinear mapping $f:M\times
M\rightarrow N$ is called \textit{non-degenerate} in both variables if $f(x,M)=0$ or $f(M,x)=0$ implies $x=0$.
We call the bilinear map $f$, \textit{``onto''} if $N$ is generated by $f(x,y)$, $x,y\in M$. We associate two
many sorted structures to every bilinear mapping described above. One of them
$$\mathfrak{U}_R(f)=\langle R,M,N,\delta,s_M,s_N\rangle,$$ where
the predicate $\delta$ describes the mapping and $s_M$ and $s_N$ describe the actions of $R$ on the modules $M$
and $N$ respectively. The other one,
$$\mathfrak{U}(f)=\langle M,N, \delta\rangle,$$
 contains only a
predicate $\delta$ describing the mapping $f$. It can be easily seen that the structure $\mathfrak{U}(f)$ is
absolutely interpretable in $\mathfrak{U}_R(f)$. We intend to show that there is a ring $P(f)$ such that
$\mathfrak{U}_{P(f)}(f)$ is absolutely interpretable in $\mathfrak{U}(f)$. Moreover this ring is the maximal
ring relative to which $f$ remains bilinear.

\subsubsection{\protect Enrichments of bilinear mappings}
In this subsection are the modules are considered to be exact. Let $M$ be an $R$-module and let $\mu:R\rightarrow P$ be an inclusion of rings. Then the $P$-module $M$ is an
\textit{$P$-enrichment} of the $R$-module $M$ with respect to $\mu$ if for every $r\in R$ and $m \in M$,
$rm=\mu(r)m$. Let us denote the set of all $R$ endomorphisms of the $R$-module $M$ by $End_R(M)$. Suppose the
$R$-module $M$ admits a $P$-enrichment with respect to the inclusion of rings $\mu:R\rightarrow P$. Then every
$\alpha\in P$ induces an $R$-endomorphism, $\phi_{\alpha}:M\rightarrow M$ of modules defined by
$\phi_{\alpha}(m)=\alpha m$ for $m \in M$. This in turn induces an injection $\phi_P:P\rightarrow End_R(M)$ of
rings. Thus we associate a subring of the ring $End_R(M)$ to every ring $P$ with respect to which there is an
enrichment of the $R$-module M.
\begin{defn}Let $f:M\times M\rightarrow N$ be an $R$-bilinear
``onto'' mapping and $\mu:R\rightarrow P$ be an inclusion of rings. The mapping $f$ admits $P$-enrichment with
respect to $\mu$ if the $R$-modules $M$ and $N$ admit $P$ enrichments with respect to $\mu$ and $f$ remains
bilinear with respect to $P$. We denote such an enrichment by $E(f,P)$.\end{defn} We define an ordering $\leq$
on the set of enrichments of $f$ by letting $E(f,P_1)\leq E(f,P_2)$ if and only if $f$ as an $P_1$ bilinear
mapping admits a $P_2$ enrichment with respect to inclusion of rings $P_1\rightarrow P_2$. The largest
enrichment $E_H(f,P(f))$ is defined in the obvious way. We shall prove existence of such an enrichment for a
large class bilinear mappings.
\begin{prop} \label{P(f)}If $f:M\times M \rightarrow N$ is a non-degenerate
``onto'' $R$-bilinear mapping over a commutative ring $R$, $f$ admits the largest enrichment.\end{prop}

\subsubsection{\protect Interpretability of the $P(f)$ structure}

Let $f:M\times M\rightarrow N$ be a non-degenerate ``onto'' $R$-bilinear mapping for some commutative ring $R$.
The mapping $f$ is said to have \textit{finite width} if there is a natural number $S$ such that for every $u\in
N$ there are $x_i$ and $y_i$ in $M$ we have
$$u=\sum_{i=1}^nf(x_i,y_i).$$
The least such number, w(f), is the \textit{width} of $f$.

A set $E=\{e_1,\ldots e_n\}$ is a complete system for $f$ if $f(x,E)=f(E,x)=0$ for $x\in M$ implies $x=0$. The
cardinality of a complete system with minimal cardinality is denoted by $c(f)$.

\textit{Type} of a bilinear mapping $f$, denoted by $\tau(f)$, is the pair $(w(f),c(f))$. The mapping $f$ is
said to be of finite type if $c(f)$ and $w(f)$ are both finite numbers.

Now we state the main theorem of this subsection:
\begin{thm}\label{ringinter}Let $f:M\times M\rightarrow $ be non-degenerate ``onto'' bilinear mapping of finite type. Then the structure $\mathfrak{U}_{P(f)}(f)$
 is absolutely interpretable in $\mathfrak{U}(f)$\end{thm}

For proofs and details we refer the reader to~\cite{alexei86} and~\cite{alexei91} although we shall describe the ring $P(f)$ later.
\section{\protect $N_{2,n}$ groups}\label{n2ns}
\subsection{\protect Definition of $N_{2,n}$ groups}

Let $R$ be a ring with unit and for an arbitrary natural number $n \geq 2$ consider the set of all
$n+(^n_2)$-tuples $((\al_i)_{1\leq i\leq n},(\ga_{ij})_{1 \leq i < j \leq n})$ of elements of $R$ where $(i,j)$ are ordered lexicographically and the same order is assumed on the $\ga_{ij}$. By
$((-),(-))$ is meant a concatenation of two tuples. We denote this set by $X$. We drop the subscripts
and denote the tuple only by $((\al_i),(\ga_{ij}))$. Always $(\bar{0})$ means that all the coordinates are $0$.
Define a multiplication on this set by:
\begin{equation}\label{multn}((\al_i),(\ga_{ij}))((\be_i),(\ga'_{ij}))=_{df}((\al_i+\be_i),(\ga_{ij}+\ga'_{ij}+\al_i\be_j)), \quad \al_i,\be_i,\ga_{ij},\ga'_{ij}\in R.\end{equation}
By $N_{2,n}(R)$ we mean the set $X$ together with the operation above.
\begin{lem}The operation defined in \eqref{multn} makes $N_{2,n}(R)$ into a group. \end{lem}
\begin{proof}
Let $x=((\al_i),(\ga_{ij}))$, $y=((\be_i),(\ga'_{ij}))$ and $z=((\de_i),(\ga''_{ij}))$ be elements of $N_{2,n}(R)$. Then
\begin{equation*}\begin{split}
(xy)z&= ((\al_i+\be_i),(\ga_{ij}+\ga'_{ij}+\al_i\be_j))z\\
&=(((\al_i+\be_i)+\de_i),(((\ga_{ij}+\ga'_{ij})+\ga''_{ij}+\al_i\be_j+(\al_i+\be_i)\de_j))\\
&=((\al_i+(\be_i+\de_i)),((\ga_{ij}+(\ga'_{ij})+\ga''_{ij})+\al_i(\be_j+\de_j)+\be_i\de_j)\\
&=((\al_i),(\ga_{ij}))((\be_i+\de_i),(\ga'_{ij}+\ga''_{ij}+\be_i\de_j))\\
&=x(yz),\end{split}\end{equation*} which proves the associativity of the operation. The identity element is
clearly $((\bar{0}),(\bar{0}))$ and if $x$ is as above then $x^{-1}=((-a_i),(a_ia_j-d_{ij}))$. So $N_{2,n}(R)$
is a group.\end{proof}
 An isomorphic copy of $N_{2,n}(R)$ is called an \textit{$N_{2,n}$ group over $R$}. If $\mathcal{R}$ is a class of
 rings,
 $N_{2,n}(\mathcal{R})$ is the class of all groups $G$ such $G\cong N_{2,n}(R)$ for some ring $R$ in $\mathcal{R}$.
 If $\mathcal{R}$ is the class of all rings a member of the class $N_{2,n}(\mathcal{R})$ is called an $N_{2,n}$ group. We note that $N_{2,2}(R)\cong UT_3(R)$.

Next proposition shows our main interest in $N_{2,n}$ groups. We postpone the proof to the end of
Subsection~\ref{Gi}.

\begin{prop}\label{freenil}If $\mathbb{Z}$ is the ring of integers then $N_{2,n}(\mathbb{Z})$ is a free 2-nilpotent group of rank
$n$.\end{prop}
\subsection{\protect Commutator subgroup and center of a $N_{2,n}$ group}

  Let $G$ be a $N_{2,n}$ group over some ring $R$ with unit. An elementary computation shows if $x=((\al_i),(\ga_{ij}))$ and $y=((\be_i),(\ga'_{ij}))$ then we have
$$[x,y]=((\bar{0}),(\al_i\be_j-\be_i\al_j)).$$ Now we can study the relation between the commutator
subgroup $[G,G]$ and the center $Z(G)$ of $G$.
\begin{lem}\label{comcenn} Let $G$ be a $N_{2,n}$ group. Then $Z(G)=[G,G]$.\end{lem}
\begin{proof}By the equation for commutators obtained above it is clear that $[G,G]$ is the set of elements of the form
$x=((\bar{0}),(\ga_{ij}))$, $\ga_{ij}\in R$, $1\leq i < j \leq n$. It is clear that for such $x$, $[x,y]=1$ for every $y\in G$. So $[G,G]\subseteq Z(G)$.
For the converse let $x=((\al_i),(\ga_{ij}))\in Z(G)$. If $y=((\be_i),(\ga'_{ij}))$ is an arbitrary element of $G$ then
we must have $[x,y]=((\bar{0}),(\al_i\be_j-\be_i\al_j))=1=((\bar{0}),(\bar{0}))$. Since this equality holds for all
elements $\be_i$ and $\be_j$ of $R$ it also holds if $\be_j=1$ and $\be_i=0$, for each $1\leq i < j \leq n$. So all
$\al_i=0$, $1\leq i \leq n-1$. Setting $\be_{n-1}=1$ and $\be_n=0$ will obtain that $\al_n=0$. So $x\in [G,G]$.
\end{proof} We note that as a consequence of Lemma~\ref{comcenn} a $N_{2,n}$ group is 2-nilpotent.

\subsection{\protect Standard basis for a $N_{2,n}$ group}
\label{generator}

When all coordinates of an element $x=((\al_i),(\ga_{ij}))$ of $N_{2,n}(R)$ are zero except possibly the $i$-th
coordinate then $x$ is denoted by $g_i^{\al_i}$. If every coordinate of $x$ is zero except possibly the $ij$-th
coordinate $x$ is denoted by $g_{ij}^{\ga_{ij}}$. In particular $g_i^0=g_{ij}^0=1$. We also assume that
$g_i^1=g_i$ for all $1\leq i \leq n$ and $g_{ij}^1=g_{ij}$ for all $1\leq i< j \leq n$. By what has been shown
above:
 $$[g_i^\al,g_j^\be]=g_{ij}^{\al\be},\quad \al,\be\in R,$$ and $[g_i,g_j]=g_{ij}$. So
$[g_i^{\al}, g_j^\be]=g_{ij}]^{\al\be}=[g_i,g_j]^{\al\be}$. Thus given an element $x=((\al_i),(\ga_{ij}))\in N_{2,n}(R)$ it is clear that
\begin{equation}\label{eq2}x=g_n^{\al_n}g_{n-1}^{\al_{n-1}}\ldots g_1^{\al_1}[g_1,g_2]^{\ga_{12}}\ldots [g_1,g_n]^{\ga_{1n}}
 [g_2,g_3]^{\ga_{23}}\ldots [g_{n-1},g_n]^{\ga_{n-1,n}}.\end{equation}
 or equivalently
 \begin{equation}\label{eq1}
 x=g_n^{\al_n}g_{n-1}^{\al_{n-1}}\ldots g_1^{\al_1}g_{12}^{\de_{12}}\ldots g_{1n}^{\de_{1n}}
 g_{23}^{\de_{23}}\ldots g_{n-1,n}^{\de_{n-1,n}}.\end{equation}

 Thus the set $\{g_i^\al|1\leq i \leq n, \al \in R\}$ is a generating set for $N_{2,n}(R)$. Moreover it should be
 clear that every element $x$ of $N_{2,n}(R)$ has a unique representation of the form given in the equations~(\ref{eq1}) and
 (\ref{eq2}). We call the elements $g_1$,\ldots, $g_n$ a \textit{standard
 basis} for the group $N_{2,n}(R)$.

 \subsection{\protect Centralizers of elements of the standard basis of a $N_{2,n}$ group} \label{Gi}

Consider a ring $R$ with unit and a $N_{2,n}$ group $G$ over $R$.
 Let $C_G(x)$ denote the centralizer of an element $x$ of $G$ in $G$.
 Now Set:
 $$G_i=_{df}C_G(g_i), \quad 1\leq i\leq n,$$
  where $g_1$,\ldots, $g_n$ constitute the standard basis for $G$.
 Let $g_i^R=\{g_i^\alpha: \alpha\in R\}$. We first prove that $G_i=g_i^R\oplus Z(G)$.
 \begin{lem}\label{formGi}For each $1\leq i\leq n$, $G_i=g_i^R\oplus Z(G)$.\end{lem}
 \begin{proof}
 Firstly we observe that $g_i^R$ is a subgroup. Let $x=((\al_{i}),(\ga_{ij}))$ be an
 arbitrary element of $G$. By the discussion in Subsection~\ref{generator}, $x=g_n\ldots g_1v$ when $v$ is an element of the center $Z(G)$.
 Then
 \begin{equation*}\begin{split}
 [g_i, x]&=[g_i, g_n^{\al_n}\ldots g_1^{a_1}v]\\
 &=[g_i,g^{\al_n}_n\ldots g^{\al_1}][g_i,v]\\
 &=[g_i,g^{\al_n}_n\ldots g^{\al_1}_1]\\
 &=g_{1i}^{-\al_i}\ldots g_{i-1,i}^{-\al_{i-1,i}}g_{i,i+1}^{\al_{i,i+1}}\ldots g_{jn}^{\al_n}
 \end{split}\end{equation*}
 For $x$ to be in $G_i$ it is necessary that $[g_i, x]=1$. So by the above equality $\al_j=0$ for $i\neq j$. So
 $x$ must have the form $g_i^\al v$ for some $a$ in $R$ and $v$ in $Z(G)$. It is also clear that $g_i^R\cap Z(G)=1$.\end{proof}
 \begin{cor}\label{Giabelian}Each $G_i$, $1\leq i \leq n$ is abelian.\end{cor}
 \begin{proof}Clear by Lemma~\ref{formGi}.\end{proof}

  Next define subgroups $G_{ij}$ of $G$ by
   \begin{equation}\label{Gij}G_{ij}=_{df}[G_i,G_j],\quad  1\leq i<j\leq n.\end{equation}
\begin{lem}\label{[Gij]} The equalities:
$$G_{ij}=[G_i,g_j]=[g_i,G_j],\quad 1\leq i < j \leq n,$$
 hold, when $G_{ij}$ are defined in Equation~\eqref{Gij}.\end{lem}
   \begin{proof}We shall prove $[g_i,G_j]=[G_i,G_j]$ for each $1\leq i<j\leq n$. The direction $\subseteq$ is obvious. For
 the converse let $x\in G_i$ and $y\in G_j$. By Lemma~\ref{formGi}, $x=g_i^av$ and $y=g_j^bv'$ for some $a,b\in
 R$ and $v,v'\in Z(G)$. Thus,
 \begin{equation*}\begin{split}
 [g_i^\al v,g_j^\be v']&=[g_i^\al,g_j^\be]\\
 &=[g_i,g_j]^{\al \be}\\
 &=[g_i,g_j^{\al \be}]\in [g_i, G_j].\end{split}
  \end{equation*}The other equality can be proved similarly. \end{proof}
 The next thing we can verify is that
\begin{equation}\label{intersectioneq}G_i\cap G_j=Z(G),\quad 1\leq i,j \leq n.\end{equation}
\begin{lem}\label{intersection} Equations~\eqref{intersectioneq} hold in for the subgroups $G_i$, $G_j$ and $Z(G)$ for each $1\leq i,j\leq n$.\end{lem}
\begin{proof}Let the element $x$ of $G$ be such that $x\in G_i$ and
$y\in G_j$. By Lemma~\ref{formGi}, $x=g_i^\al v=g_j^\be v'$ for some $\al, \be\in R$ and $v,v'\in Z(G)$, implying that
$\al=\be=0$. Therefore $x\in Z(G)$. The other direction is clear.\end{proof} We assemble the lemmas and corollaries
above in a single proposition.
\begin{prop}\label{basen2}Let $G$ be a $N_{2,n}$ group over a ring $R$ with unit. Suppose $g_1$,\ldots,$g_n$ constitute the standard
basis for $G$. Let $G_i=C_G(g_i)$, $1\leq i \leq n$, and $G_{ij}=[G_i,G_j]$, $1\leq i < j \leq n$.
 Then the
following statements hold,
\begin{enumerate}
\item $G_{ij}=[G_i,g_j]=[g_i,G_j],\quad 1\leq i<j\leq n$, \item $G_i \cap G_j=Z(G), \quad 1\leq i<j\leq n$,
\item $[G_i,G_i]=1, \quad 1\leq i\leq n$,
  \item$[G,G]=Z(G)$,
  \item every element of $G$ can be written as $u_n u_{n-1}\ldots u_1v$ where each $u_i \in G_i$
and $v\in Z(G)$ and each $u_i$ is unique modulo the center. Moreover each $v \in Z(G)$ can be uniquely written
as $u_{12}\ldots u_{1n}u_{23}\ldots u_{n-1,n}$ when $u_{ij}\in G_{ij}$.
\end{enumerate}\end{prop}
\begin{proof} See Subsection~\ref{generator}, Lemmas~\ref{[Gij]},~\ref{intersection} and Corollary~\ref{Giabelian}.\end{proof}

 Now we are in a good set up to prove Proposition~\ref{freenil}.

 \begin{proof}(\textbf{Proof of Proposition~\ref{freenil}})
  Let $F(n)$ be the free group on generators $u_1, \ldots, u_n$. Let $\Gamma_3(F(n))$ denote
the third term of the lower central series of $F(n)$. Let $g_1$, \ldots, $g_n$ be elements of the standard basis
for $N_{2,n}(\mathbb{Z})$. Note that $\{g_1, \ldots, g_n\}$ is a generating set for $N_{2,n}(\mathbb{Z})$. The
mapping:
$$F(n)/\Gamma_3(F(n)) \longrightarrow N_{2,n}(\mathbb{Z}), \quad u_i\Gamma_3(F(n)) \mapsto g_i, \quad 1\leq i \leq n,$$
is a well defined homomorphism of groups since $\Gamma_3(F(n))$ is generated by the simple commutators
$[[u_i,u_j],u_k]$, $i\neq j$, and $[[g_i,g_j],g_k]=1$ holds in $N_{2,n}(\mathbb{Z})$. It is also a
surjection since $g_1$,\ldots, $g_n$ generate $N_{2,n}(\mathbb{Z})$.

We prove that it is also an injection. Notice that for every integer $k$ and $m$,
\begin{equation}\label{elim}[u_i^m\Gamma_3(F(n)),u_j^k\Gamma_3(F(n))]=[u_i,u_j]^{mk}\Gamma_3(F(n))\end{equation} holds for each
pair $u_i$, $u_j$ of elements in $\{u_1,\ldots,u_n\}$. So every element $u$ in $F(n)/\Gamma_3(F(n))$ can be
brought to the form in (\ref{eq2}), $g_i$ substituted by $u_i$ for each $1\leq i \leq n$. This can be done using
the so-called collection process, applying the relations~\eqref{elim} and using the fact that all the
commutators in $F(n)/\Gamma_3(F(n))$ belong to the center of the group. For example if $u_i$ and $u_j$
are such that $i<j$ and $m$ and $k$ are integers then:
\begin{equation*}\begin{split}
u_i^mu_j^k\Gamma_3(F(n))&=u_i^ku_j^m[u_i^m,u_j^k]\Gamma_3(F(n))\\
&=u_j^ku_i^m[u_i,u_j]^{mk}\Gamma_3(F(n)).\end{split}\end{equation*} By repeating this process
finitely many times and moving the commutators to the right hand side we arrive at the indicated form for any
element of $F(n)/\Gamma_3(F(n))$. Therefore under the mapping defined above $u$ gets
mapped to an element $g$ of $N_{2,n}(\mathbb{Z})$ with a representation exactly like what appears in
(\ref{eq2}). This form is unique so the element $g$ is trivial in $N_{2,n}(\mathbb{Z})$ if and only if all the
exponents in (\ref{eq2}) are zero if and only if $u$ is trivial in $F(n)/\Gamma_3(F(n))$. And we are done.
\end{proof}
\begin{rem}\label{Rgr}Assume $\{g_1, \ldots, g_n\}$ is the standard basis for $G=N_{2,n}(\mbb{Z})$. Then it is not too hard to verify that the tuple $(g_n, g_{n-1}, \ldots , g_1, g_{12}, \ldots , g_{n,n-1})$ is a Mal'cev basis for $G$. So if $R$ is a binomial domain by definition $N_{2,n}(R)$ is the Mal'cev $R$-completion of $G$ an hence an $R$-group. We'll use this fact in the last section.\end{rem}
\section{\protect $QN_{n,2}$ groups}\label{quasiqnn2}
\subsection{\protect Definition of $QN_{n,2}$ groups}\label{quasigroup}

Let $f \in S^2(R^+, \oplus_{i=1}^{(^n_2)}R^+$ be a symmetric 2-cocycle, where $R$ is a ring with unit. Such a
2-cocycle has $(^n_2)$ coordinates $f_{jk}:R^+\times R^+ \longrightarrow R$, $1\leq j < k \leq n$, each of which
is a symmetric 2-cocycle. Now for each $i$, $1\leq i\leq n$, let $f^i:S^2(R^+, \oplus_{i=1}^{(^n_2)}R^+$
be a symmetric 2-cocycle with components $f^i_{jk}$, $1\leq j < k \leq n$.

We define a new multiplication $\odot$ on the underlying set $X$ of $N_{2,n}(R)$ by
\begin{equation}\label{multqn}((\al_i),(\ga_{ij}))\odot((\be_i),(\ga'_{ij}))=((\al_i+\be_i),(\ga_{ij}+\ga'_{ij}+\al_i\be_j+\sum_{k=1}^nf^k_{ij}(\al_k,\be_k))).\end{equation}
\begin{lem}The set $X$ is a group with respect to the multiplication $\odot$ defined in~\eqref{multqn}.\end{lem}
\begin{proof}
Let $g^i:\oplus_{i=1}^nR^+\times \oplus_{i=1}^nR^+ \rightarrow \oplus_{i=1}^{\binom{n}{2}}$ be defined by
$$g^i((\alpha_1, \ldots , \al_n), (\beta_1 ,\ldots , \beta_n))=f^i(\alpha_i,\beta_i).$$
It is easy to verify that $g^i\in B^2(\oplus_{i=1}^nR^+, \oplus_{i=1}^{\binom{n}{2}})$. Now it is clear that the $(X,\odot)$ is a central extension of $\oplus_{i=1}^nR^+$ by $\oplus_{i=1}^{\binom{n}{2}}R^+$ via a 2-cocycle $f=(f_{ij})_{1\leq i < j \leq n}$ defined by
$$f_{ij}((\alpha_1, \ldots , \al_n), (\beta_1 ,\ldots , \beta_n))=\alpha_{ij}+\sum_{i=1}^ng^i_{ij}((\alpha_1, \ldots , \al_n), (\beta_1 ,\ldots , \beta_n)).$$\end{proof}

We denote the new group by $N_{2,n}(R,f^1 \ldots
f^n)$. If $\mathcal{R}$ is a class of rings with unit, by $QN_{2,n}(\mathcal{R})$ we mean the class of all
groups $G$ such that $G\cong N_{2,n}(R,f^1,\ldots,f_n)$ for some ring $R$ in $\mathcal{R}$ and symmetric
2-cocycles $f^i:R^+\times R^+\rightarrow R^{(^n_2)}$, $i=1,\ldots,n$. Such a group $G$ is called a $QN_{2,n}$
group over $R$. If $\mathcal{R}$ is the class of all rings a member of the class $QN_{2,n}(\mathcal{R})$ is
called a $QN_{2,n}$ group.
\subsection{\protect Commutator subgroup and center of a $QN_{2,n}$ group}

Let $G$ be a $QN_{2,n}$ group over a ring $R$ with unit. To give a formula for the commutator of two elements we
need  to verify a basic fact about symmetric 2-cocycles.

Let $x=((\al_i),(\ga_{ij})$ and $y=((\be_i),(\ga'_{ij}))$ be in $G$ then
\begin{equation*}\begin{split}
x^{(-1)}\odot y^{(-1)}\odot x \odot y&= ((\bar{0}),(\ga_{ij}+\ga'_{ij}-\ga_{ij}-\ga'_{ij}\\
&\quad +\al_i\al_j+\be_i\be_j+\al_i\be_j + \al_i\be_j + (-\al_i-\be_i)(\al_j+\be_j)\\
&\quad - \sum_{k=i}^n(f^k_{ij}(\al_k,-\al_k) -f^k_{ij}(\be_k,-\be_k)\\
&\quad + f^k_{ij}(\al_k,\be_k)+f^k_{ij}(-\al_k-\be_k,\al_k+\be_k))\\
&= ((\bar{0}),(\al_i\be_j-\be_i\al_j)) \end{split}\end{equation*} by the above lemma. Thus commutators in $QN_{2,n}$ and
$N_{2,n}$ groups coincide. So we have the lemma:
\begin{lem}\label{comcenq}In a $QN_{2,n}$ group $G$, $Z(G)=[G,G]$.\end{lem}
\begin{proof}The proof goes through exactly like that of Lemma~\ref{comcenn}.\end{proof}
\subsection{\protect Standard basis for a $QN_{2,n}$ group}

Again as in the $N_{2,n}$ groups we denote an element $((\alpha_i),(\gamma_{ij}))$ which has zeros everywhere except
possibly at the $i$-th position by $g_i^{a_i}$ and the one which has zeros everywhere except possibly at $ij$-th
position by $g_{ij}^{\ga_{ij}}$. We call the set $\{g_1, \ldots ,g_n\}$ the \textit{standard basis} of the group
$QN_{2,n}(R)$. Let us note that for a $QN_{2,n}$ group $G$ over a ring $R$ with unit and the standard basis
$\{g_1,\ldots ,g_n\}$, the quotient $G/Z(G)$ is a free module over $R$ of rank $n$ generated by
$\{g_1Z(G),\ldots ,g_nZ(G)\}$. Moreover $Z(G)=[G,G]$ is a free $R$-module of rank $(^n_2)$ generated
by the $g_{ij}=[g_i,g_j]$, $1\leq i < j \leq n$.

\begin{prop}\label{baseqn2}Let $G$ be a $QN_{2,n}$ group over a ring $R$ with unit, $G_i$ for each $1\leq i\leq n$ and $G_{ij}$
  for each $1\leq i<j \leq n$ be defined as in proposition~\ref{basen2}. Then all the conditions (1)-(5) in
  proposition~\ref{basen2} are also true in the group $G$.\end{prop}
  \begin{proof}Similar to the proof of Proposition~\ref{basen2}.\end{proof}
  \subsection{\protect Generators and relations for a $QN_{2,n}$ group}

  Here we specify a set of generators and relations for a $QN_{2,n}$ group.
  \begin{lem}\label{salam}The group $G=N_{2,n}(R,f^1, \ldots ,f^n)$ is
    generated by
    $$\{g_i^a,g_{kl}^b:1\leq i\leq
n, 1\leq k<l\leq n,\quad \alpha,\beta\in R\}, $$ and
    defined by the relations:\\
    (a) $[g_i^\alpha,g_j^\beta]=g_{ij}^{\alpha\beta}$,\quad for all $1\leq i < j\leq n$, $\alpha,\beta\in R$\\
    (b) $[g_i^\alpha,g_{kl}^\beta]=[g_{rs}^\al, g_{kl}^\be]=1$, \quad for all $1\leq i \leq n$, $1\leq k < l \leq n$, $1\leq r < s \leq n$ and $\alpha,\beta\in R$,\\
    (c) $g_i^\alpha \odot g_i^\beta=g_i^{(\alpha+\beta)}g_{12}^{f^i_{12}(\alpha,\beta)} \ldots
    g_{n-1,n}^{f^i_{n-1,n}(\alpha,\beta)},$\quad for all $1\leq i\leq n$,
    $\alpha,\beta\in R$\\
    (d) $g_{ij}^\alpha\odot g_{ij}^\beta=g_{ij}^{\alpha+\beta}$,\quad for all $1\leq i < j \leq n$, $\alpha,\beta\in R$.
    \end{lem}
\begin{proof} clearly the set
$$\mathcal{G}=\{g_i^a,g_{kl}^b|1\leq i\leq
n, 1\leq k<l\leq n,\quad\alpha,\beta\in R \},$$ is a generating set for $G$. Let $F$ be the free group generated
by the set $\mathcal{G}$ and $\mathcal{R}$ be the normal subgroup of $F$ generated by the relations (a)-(d)
above, multiplication $\odot$ taken to be concatenation. Now consider the group $\langle \mathcal{G}|\mathcal{R}
\rangle$, the quotient of $F$ by $\mathcal{R}$. Consider the mapping:
$$\langle \mathcal{G}|\mathcal{R}\rangle \longrightarrow G, \quad g^\alpha_i \mapsto g_i^\alpha,
\quad g_{kl}^\alpha \mapsto g_{kl}^\alpha$$ for every $\alpha \in R$, $1\leq i \leq n$ and $1\leq k <l \leq n$.
The map is a well-defined homomorphism since all the relations (a)-(d) hold in $G$. Every word $W$ in $\langle
\mathcal{G}|\mathcal{R}\rangle$ is equivalent to a word with the form given in (\ref{eq1}), multiplication taken
to be concatenation. this element gets mapped to an element $g$ with the same form in the group $G$,
multiplication taken to be $\odot$. This form is unique by Proposition~\ref{baseqn2}. So $g$ is trivial in $G$
if and only if $W$ is trivial in $\langle \mathcal{G}|\mathcal{R}\rangle$. The proposition is proved.
\end{proof}

\section{\protect $QN_{2,n}$ groups over commutative rings}\label{commutsection}
\label{rbil}

Let $G$ be a nilpotent group of class 2. We associate to $G$ a bilinear map
$$f_G:\ G/[G,G] \times G/[G,G] \longrightarrow [G,G], \quad (x[G,G],y[G,G])\mapsto [x,y]$$
$x,y \ in G$. Note that in a $QN_{2,n}$ group $G$, $[G,G]=Z(G)$. So the bilinear map $f_G$ becomes
$$f_G: G/Z(G)\times G/Z(G), \quad (xZ(G),yZ(G))\mapsto [x,y]$$
$x,y\in G$.

When $G$ is a $QN_{2,n}$ the bilinear map $f_G$ is both ``onto'' and non-degenerate so as in Section~\ref{modelbilin}we can associate to it a commutative associative ring $P(f_G)$ with unit relative to which $f_G$ is bilinear and $P(f_G)$ is the maximal such ring. Actually $P(f_G)$ is the set of all pairs
$$(\phi_1,\phi_0)\in E= End(G/Z(G))\times End(Z(G)),$$ where $End(-)$ refers to the endomorphism ring of the corresponding group, which satisfy the identity
$$f_G(\phi_1(xZ(G)),yZ(G))=f_G(xZ(G),\phi_1(yZ(G)))=\phi_0(f_G(x,y))$$
for all $x,y \in G$.
\begin{lem}Let $G$ be a $QN_{2,n}$ group over a commutative associative ring $R$
with unit. Let $(\phi_1,\phi_0)\in P(f_G)$ and $x,y \in G$. Then $\phi_1(xZ(G))=(xZ(G))^\gamma $ and
$\phi_0(f_G(x,y))=f_G(x,y)^\gamma$ for some $\gamma \in
R$.\label{action}\end{lem}
\begin{proof} Let $\{g_1,g_2,\ldots ,g_n \}$ be the standard basis for $G$ and $g_{ij}=[g_i,g_j]$  for  $1\leq  i<j\leq n$.
 Since $g_iZ(G)$, $1\leq i\leq n$, generate $G/Z(G)$ and $g_{ij}$, $1 \leq i<j \leq n$ generate $Z(G)$ it is enough to study the action of the $\phi_i$ on powers of basis elements.

Fix $1\leq i \leq n$. We assume that $\phi_1(g_iZ(G))=g_1^{\alpha_1}g_2^{\alpha_2}\ldots g_n^{\alpha_n}Z(G)$
 for some $\alpha_k \in R$, $1\leq k \leq n$.  Now
\begin{equation*}\begin{split}
g_{1j}^{\alpha_1}\ldots g_{j-1,j}^{\alpha_{j-1}}g_{j,j+1}^{-\alpha_{j+1}}\ldots g_{jn}^{-\alpha_n} &=
[g_1^{\alpha_1}\cdots
g_n^{\alpha_n},g_i]\\
&= f_G(\phi_1(g_iZ(G)),g_iZ(G))\\
&=\phi_0([g_i, g_i])\\
&= 1
\end{split}\end{equation*}
So $\alpha_k=0$, $k\neq i$. Set $\alpha_k=\ga$. Now pick $1\leq j \leq n$ such that $i\neq j$. By an argument like the one above we may conclude that there exists $\de\in R$  such that $\phi_1(g_jZ(G))=(g_{j}Z(G))^{\de}$. Then
\begin{align*}
 g_{ij}^\ga &= f_G(\phi_1(g_i), g_j)\\
 &= f_G(g_i, \phi_1(g_j))\\
 &= g_{ij}^\de
\end{align*}
So $\ga=\de$. This proves that for any $1\leq k \leq n$, $\phi_1(g_kZ(G))=(g_kZ(G))^\ga$.

Next let $\alpha \in R$ and pick $1 \leq i < j \leq n$. Suppose $$\phi_1(g_i^\alpha Z(G))=g_1^{\alpha_1}\ldots
g_n^{\alpha_n}Z(G).$$ Then,
\begin{equation*}\begin{split}
g_{1j}^{\alpha_1}\cdots g_{j-1,j}^{\alpha_{j-1}}g_{j,j+1}^{-\alpha_{j+1}}\ldots g_{jn}^{-\alpha_n} &=
[g_1^{\alpha_1}\cdots g_n^{\alpha_n},g_j] \\
&=f_G(\phi_1(g_i^\alpha Z(G)),g_jZ(G))\\
& =f_G((g_iZ(G))^\alpha,\phi_1(g_jZ(G))) \\
&=[g_i^\alpha,g_j^\gamma]=g_{ij}^{\alpha\gamma} \end{split}\end{equation*} So $\alpha_i=\alpha\gamma$ and
$\alpha_k=0$ if $k\neq i,j$. To prove $\alpha_j=0$ it is enough to consider $f_G(\phi_1((g_iZ(G))^{\alpha}),g_iZ(G))=1$. It is also clear that $\phi_0(g_{ij}^\alpha)=g_{ij}^{\al\ga}$ for all $1\leq i < j \leq n$.

Since $\phi_1$ and $\phi_0$ are endomorphisms the statement is clear now.
\end{proof}
\begin{prop}\label{P(f)-R}Let $R$ be a commutative associative ring with unit and $G$ be a $QN_{2,n}$ group over $R$. Then
$P(f_G)\cong R$\end{prop}
\begin{proof}Define a mapping
$$\eta:P(f_G) \rightarrow R, \quad (\phi_1,\phi_0)\mapsto \gamma_{\phi}$$
where $\phi_1(xZ(G))=(xZ(G))^{\gamma_{\phi}}$ for $x\in G$ and $\phi_0(y)=y^{\gamma_{\phi}}$ for $y\in
Z(G)$. Such a $\gamma_{\phi}$ exists by Lemma~\ref{action}. The mapping is well defined since if
$(xZ(G))^{\gamma_1}=(xZ(G))^{\gamma_2}$ for all $x\in G$ then also $(g_iZ(G))^{\gamma_1}=(g_iZ(G))^{\gamma_2}$ which implies
$\gamma_1=\gamma_2$.
Let $\gamma$ be an element of the ring $R$. Define a triple $(\phi_1,\phi_1,\phi_0)$ where $\phi_1 \in
End(G/Z(G))$ and $\phi_0 \in End(Z(G))$ by setting $\phi_1(xZ(G))=(xZ(G))^{\gamma}$, for $x\in G$ and
$\phi_0(z)=z^{\gamma}$, for $z\in Z(G)$. We show that $(\phi_1,\phi_1,\phi_0)\in P(f_G)$. Let $ \{g_1,\ldots,
g_n \}$ be the standard basis for $G$, $g_{ij}=[g_i,g_j]$, $1\leq i<j \leq n$, $xZ(G)=(g_1Z(G))^{\alpha_1}\cdots
(g_nZ(G))^{\alpha_n}$ and $yZ(G)=(g_1Z(G))^{\beta_1}\cdots (g_nZ(G))^{\beta_n}$ for $x,y\in G$. Then by associativity and
commutativity of $R$,
\begin{equation*}\begin{split}
f_G(x^{\gamma}Z(G),yZ(G))&=f_G(g_1^{\alpha_1}\ldots g_n^{\alpha_n}Z(G),g_1^{\beta_1}\ldots g_n^{\beta_n}Z(G))\\
&=g_{12}^{((\alpha_1\gamma)\beta_2-\beta_1(\alpha_2\gamma))}\ldots
g_{n-1,n}^{((\alpha_{n-1}\gamma)\beta_n-\beta_{n-1}(\alpha_n\gamma))}\\
[x,y]^{\gamma}&=(g_{12}^{(\alpha_1\beta_2-\beta_1\alpha_2)}\ldots
g_{n-1,n}^{(\alpha_{n-1}\beta_n-\beta_{n-1}\alpha_n)})^{\gamma}\\
&= g_{12}^{(\alpha_1(\beta_2\gamma)-\beta_1(\gamma\alpha_2))}\ldots
g_{n-1,n}^{(\alpha_{n-1}(\beta_n\gamma)-(\beta_{n-1}\gamma)\alpha_n)}\\
&=f_G(xZ(G),y^{\gamma}Z(G)).
\end{split}\end{equation*}
So $(\phi_1,\phi_0)\in P(f_G)$. This proves the surjectivity of $\eta$. If $(\phi_1,\phi_0)\in
P(f_G)$ maps to the zero of the ring $R$ under the mapping $\eta$, it means that $\phi_1$  and $\phi_0$ are zero
endomorphisms. Hence $(\phi_1,\phi_0)$ is the zero of $P(f_G)$. Hence the mapping $\eta$ is injective.

To prove that $\eta$ is an additive homomorphism note that
$$\eta((\phi_1,\phi_0)+(\psi_1,\psi_0))=\eta((\phi_1+\psi_1,\phi_0+\psi_0))$$
is an element $\gamma$ of $R$ such that for every $x\in G$,
$$(xZ(G))^\gamma=\phi_1(xZ(G))\psi_1(xZ(G))$$ and for every $y\in Z(G)$, $y^{\gamma}=\phi_0(y)\psi_0(y)$. But
$$\phi_1(xZ(G))\psi_1(x)=(xZ(G))^{\gamma_{\phi}}(xZ(G))^{\gamma_{\psi}}$$ and
$\phi_0(y)\psi_0(y)=y^{\gamma_{\phi}}y^{\gamma_{\psi}}$. Thus $\gamma=\gamma_{\phi}+\gamma_{\psi}$, hence $\eta$
is an additive homomorphism. On the other hand identities
$\phi_1\circ\psi_1(xZ(G))=(xZ(G))^{\gamma_{\psi}\gamma_{\phi}}=(xZ(G))^{\gamma_{\phi}\gamma_{\psi}}$ and
$\phi_0\circ\psi_0(y)=y^{\gamma_{\psi}\gamma_{\phi}}=y^{\gamma_{\phi}\gamma_{\psi}}$ imply that $\eta$ is a multiplicative homomorphism. The proposition is proved.\end{proof}

\begin{prop} \label{ringiso10}Let $\varphi: G \rightarrow H$ be an isomorphism of $QN_{2,n}$ groups. Then $P(f_G)\cong P(f_H)$.\end{prop}
\begin{proof}The isomorphism $\varphi$ induces an isomorphism $\varphi_1:G/Z(G)\rightarrow H/Z(H)$ and restricts to an isomorphism $\varphi_0:Z(G)\rightarrow Z(H)$. We claim that if $(\phi_1,\phi_0)\in P(f_G)$ then
$$(\varphi_1\phi_1\varphi_1^{-1}, \varphi_0\phi_0\varphi_0^{-1})\in P(f_H)$$ and
$$\theta: P(f_G)\rightarrow P(f_H), \quad (\phi_1,\phi_0)\mapsto (\varphi_1\phi_1\varphi_1^{-1}, \varphi_0\phi_0\varphi_0^{-1}),$$
is an isomorphism of rings.\\
Let $(\phi_1,\phi_0)\in P(f_G)$, $x,y \in G$, $\varphi(x)=z$ and $\varphi(y)=t$. We note that
\begin{equation*}\begin{split}
f_H(\varphi_1(xZ(G)), \varphi_1(yZ(G))&=f_H(zZ(H), tZ(H))\\
&=[z,t]=[\varphi(x),\varphi(y)]\\
&=\varphi([x,y])=\varphi_0([x,y])\\
&=\varphi_0(f_G(xZ(G),yZ(G)))\end{split}\end{equation*}
Next using this equality we get
\begin{equation*}\begin{split}
f_H(\varphi_1\phi_1\varphi_1^{-1}(zZ(H)), tZ(H))&= f_H(\varphi_1\phi_1(xZ(G)),\varphi_1(yZ(G)))\\
&=\varphi_0(f_G(\phi_1(xZ(G),y(Z(G)))\\
&=\varphi_0\phi_0([x,y])\\
&=\varphi_0\phi_0\varphi_0^{-1}([z,t]).\end{split}\end{equation*}
With a similar argument we also get that $$f_H(zZ(H), \varphi_1\phi_1\varphi_1^{-1}( tZ(H)))= \varphi_0\phi_0\varphi_0^{-1}([z,t])$$ which proves that $(\varphi_1\phi_1\varphi_1^{-1},\varphi_0\phi_0\varphi_0^{-1}) \in P(f_H)$. \\
To prove that $\theta$ is surjective we should follow an argument similar to the above. Injectivity of $\theta$ as well as it being a ring homomorphism is clear.\end{proof}
\begin{thm}\label{commut}Let $R$ and $S$ be commutative associative rings with unit.
If $$N_{2,n}(R,f_1,\ldots f_n)\cong N_{2,n}(S,q_1, \ldots q_n)$$
as groups then $R\cong S$ as rings.\end{thm}
\begin{proof} This is a direct corollary of Lemma~\ref{P(f)-R} and Proposition~\ref{ringiso10}.\end{proof}
\section{Characterization of $Q_{2,n}$ groups over commutative associative rings}\label{charthmsec}
\subsection{ groups with a basis}

\begin{defn}[Basis] Let $H$ be a group with distinct nontrivial elements $h_1$, $h_2$,\ldots, $h_n$. Let $H_1$, $H_2$,
\ldots, $H_n$ and $H_{12}$, \ldots, $H_{n-1,n}$ be subgroups of $H$ satisfying the following conditions:
\begin{enumerate}
\item$H_i=C_H(h_i)$ and $H_{ij}=[h_i,H_j]=[H_i,h_j]=[H_i,H_j]$, $1\leq i<j \leq n$, \item$H_i \cap H_j=Z(H)$,
$1\leq i<j \leq n$, \item $[H_i,H_i]=1$, $1 \leq i \leq n$, \item$[H,H]\subseteq Z(H)$, \item \begin{itemize}
\item[(a)] every element of $H$ can be written as $u_n u_{n-1}\ldots u_1v$ where each $u_i \in H_i$ and $v\in
Z(H)$ and each $u_i$ is unique modulo the center, \item[(b)] each $v \in Z(H)$ can be uniquely written as
$u_{12}\ldots u_{1n}u_{23}\ldots u_{n-1,n}$ when $u_{ij}\in H_{ij}$. \end{itemize}\end{enumerate}Then
$\fra{b}=\{h_1, h_2, \ldots ,h_n \}$ is called a basis for $H$.\label{quasibasis}\end{defn}
\begin{defn}[$P(f_H)$-basis]Let $H$ be a group with a set of elements $\fra{b}=\{h_1,\ldots , h_n\}$. We call $\fra{b}$ a $P(f_H)$-basis for $H$ if
\begin{itemize}
\item $\fra{b}$ is a basis for $H$,
\item $H_i/Z(H)$ is a cyclic $P(f_H)$-module generated by $h_iZ(H)$, for all $1\leq i \leq n$.
\end{itemize}\end{defn}
\begin{cor}\label{cexact} Each $H_i/Z(H)$, $1\leq i \leq n$, is a torsion free $P(f_H)$-module. Moreover each $H_{ij}$ is also a torsion free $P(f_H)$ module generated by $[h_i,h_j]$.\end{cor}
\begin{proof} Pick $1\leq i \leq n$ and suppose that there is an element $\alpha\in P(f_H)$ such that $(xZ(H))\al=0$ for all $x\in H_i$. Pick $x\in H_i$ such that $x \notin Z(H)$. Let $i\neq j$, $1\leq j \leq n$ and pick any $y\in H_j$ where $y\notin Z(H)$. Then
 $$1=f_H((xZ(H))^\alpha, yZ(H))=f_H(xZ(H), (yZ(H))^\al )=[x,y'],$$
 where $y'\in H_j$ and $y'Z(H)=(yZ(H))^\al$. So $y'\in H_i\cap H_j =Z(H)$. Since $j$ was arbitrary by condition $5-a$ of the definition of a basis $\alpha \in Ann_{P(f_H)}(H/Z(H))={0}$. The first statement follows. The second statement should also be clear now.
\end{proof}

\begin{lem}\label{almost}Let $H$ be a group with elements $h_1$,\ldots,$h_n$ constituting a basis for $H$. Then the ring
$P(f_H)$ and its action on $H/Z(H)$ and $Z(H)$ are absolutely interpretable in $H$.\end{lem}
\begin{proof}The bilinear map $f_H$ has width $\frac{n(n-1)}{2}$. Moreover the set $\{h_1,\ldots h_n\}$ is a finite
complete system for $f_H$. So by Theorem ~\ref{ringinter} the structure
$$\mathfrak{U}_{P(f_H}(f_H)=\langle
P(f_H),H/Z(H), Z(H), s_{H/Z(H)},s_{Z(H)}, \delta_{f_H}\rangle,$$ where $s_{H/Z(H)}$ and $s_{Z(H)}$ describing
the action of $P(f_H)$ on $H/Z(H)$ and $Z(H)$ respectively is absolutely interpretable in
$$\mathfrak{U}(f_H)=\langle H/Z(H), Z(H), \delta_{f_H}\rangle.$$ The factor group $H/Z(H)$ is absolutely
interpretable in $H$. The subgroup $Z(H)$ is clearly definable without parameters.
There is a formula of signature of groups describing the bilinear mapping $f_H$, since $f_H$ is defined just by
commutators, $Z(H)$ is absolutely definable in $H$ and $H/Z(H)$ is absolutely interpretable in $H$. So $\mathfrak{U}(f_H)$ is absolutely interpretable in $H$. Therefore $\mathfrak{U}_{P(f_H)}(f)$ is
absolutely interpretable in $H$.\end{proof}
\begin{lem}\label{Pfhbasecor}Let $G$ be a $QN_{2,n}$ group over $R$ with a $P(f_G)$-basis $$\fra{b}=\{g_1,\ldots , g_n\}.$$ If $\varphi: G \rightarrow H$ be an isomorphism of groups then the image $\fra{b}^{\varphi}$ of $\fra{b}$ under $\varphi$ is a $P(f_H)$-basis for $H$.\end{lem}
\begin{proof} Let $\varphi(g_i)=h_i$ for each $1\leq i \leq n$. It is clear that $\fra{b}^{\varphi}$ is a basis. It remains to prove that each $C_H(h_i)/Z(H)=(h_iZ(H))^{P(f_H)}$, $1\leq i \leq n$. Let $\phi\in End(G/Z(G))$ and $\varphi_1: G/Z(G)\rightarrow H/Z(H)$ be the isomorphism induced by $\varphi$. Let $xZ(G)=\phi(g_iZ(G))$, $x\in C_G(g_i)$  and $y=\varphi(x)$. Then
\begin{equation*}\begin{split}
\varphi_1\phi\varphi^{-1}_1(h_iZ(H))&= \varphi_1\phi (g_iZ(G))\\
&=\varphi_1(xZ(G))\\
&= yZ(G).\end{split}\end{equation*}
Having this and comparing it with the isomorphism $$\theta: P(f_G)\rightarrow P(f_H)$$ defined in the  proof of Proposition~\ref{ringiso10} the result is clear. \end{proof}

\subsection{Characterization theorem}

\begin{thm}[Characterization theorem]\label{charthm}  Let $h_1,
\ldots, h_n$ be some elements of $H$. Then the following are equivalent:\\
(1) The group $H$ has a $P(f_H)$-basis, $\fra{b}=\{h_1, \ldots ,h_n\}$; \\
(2) There is a commutative associative ring $R$ and symmetric 2-cocycles
$$f^i:R^+\times R^+\rightarrow \oplus_{i=1}^{\binom{n}{2}}R^+,$$
such that $$(H, \fra{b}) \cong N^*_{2,n}(R,f^1, \ldots, f^n)=(G,\fra{b}')$$
via an isomorphism $\varphi:(G,\fra{b}') \rightarrow (H,\fra{b})$ of enriched groups where $\fra{b}'=\{g_1, \ldots , g_n \}$ is the standard basis for $G$.\\
If (1) holds each symmetric cocycle $f^i:R^+\times R^+\rightarrow \oplus_{i=1}^{\binom{n}{2}}R^+$, $1\leq i\leq n$ is constructed
such that each $H_i=C_H(h_i)$ is an abelian extension of $\oplus_{i=1}^{\binom{n}{2}}R^+$ by $R^+$ via the symmetric 2-cocycle
$f^i$.\end{thm}
\begin{proof}
(2) $\Rightarrow$ (1) By Proposition~\ref{baseqn2} $\fra{b}'$ is a basis of $G$. Since $\fra{b}'$ is the standard basis of $G$ it is clear that each $C_G(g_i)/Z(G)$ is a cyclic $R$-module. By proof of Proposition~\ref{P(f)-R}, for each element $\alpha$ of $R$ there is a unique element of $P(f_G)$  acting on $G/Z(G)$ by $\alpha$. So each $C_G(g_i)$ is a cyclic $P(f_G)$ module. So $\fra{b}'$ is a $P(f_G)$ basis for $G$. Now Lemma~\ref{Pfhbasecor} implies that $\fra{b}$ is a $P(f_H)$-basis for $H$. \\
To prove (1)$\Rightarrow$ (2) we prove that the relations
(a)-(d) of Lemma~\ref{salam} hold  with a suitable choice for $h_i^\alpha$ among the representatives
$(h_iZ(H))^\alpha$ for $\alpha\in R$ and $1\leq i\leq n$.
 Consider the following maps:
$$\tau_{i,ij}:H_i\rightarrow H_{ij}, \quad x\mapsto [x,h_{j}],$$
and
$$\tau_{j,ij}:H_j\rightarrow H_{ij}, \quad x\mapsto [h_i,x].$$
for each $1\leq i < j \leq n$. All the $\tau_{i,ij}$ and $\tau_{j,ij}$ are group homomorphisms by condition (4) of Definition~\ref{quasibasis} of a basis. They are surjective by condition (1) and they all have the same kernel, $Z(H)$ by (2) of definition of the basis. So $\tau_{i,ij}$ induces an isomorphism between $H_i/Z(H)$ and $H_{ij}$ and $\tau_{j,ij}$ induces an isomorphism between  $H_j/Z(H)$ and $H_{ij}$. So as each $H_i/Z(H)$ is a cyclic $P(f_H)$-module generated by $h_iZ(H)$, the element $h_{ij}$ generates $H_{ij}$ as a cyclic $P(f_H)$-module.
Set $R=P(f_H)$. Then by  Corollary~\ref{cexact},
$$\mu_{ij}: h_{ij}^R\rightarrow R^+, \quad h_{ij}^\zeta \mapsto \zeta,$$
 is a group isomorphism. Moreover (5)-(b) of the definition of a basis together with above explanations implies that there is an isomorphism $\eta:\oplus _{i=1}^{(^n_2)}R^+ \rightarrow Z(H)$. Now consider the following sequences of abelian groups:
$$0\rightarrow \oplus^{(^n_2)}_{i=1}R^+  \xrightarrow{ \eta} H_i \xrightarrow{ \mu_{i,i+1}\tau_{i,ii+1}} R^+ \rightarrow 0,$$
where $1\leq i \leq n-1$ and
$$0\rightarrow \oplus_{i=1}^{(^n_2)}R^+ \xrightarrow{\eta } H_n \xrightarrow{\mu_{n-1,n}\tau_{n, (n-1,n)}} R^+ \rightarrow 0.$$
 They are clearly exact. Let for each $1\leq i\leq n$, $f^i:R^+ \times
R^+\rightarrow \oplus_{i=1}^{\binom{n}{2}}R^+$ be the 2-cocycle corresponding to the extension above. Each $f^i$ is clearly a symmetric 2-cocycle since $H_i$ is abelian by condition (3)
of Definition ~\ref{quasibasis}, hence $H_i$ is an abelian extension of $\oplus_{i=1}^{(^n_2)}R^+$ by $R^+$ via the 2-cocycle
$f^i$. Therefore $H_i\cong R^+\times \oplus_{i=1}^{\binom{n}{2}}R^+=K_i$, $1\leq i\leq n$, as groups when the multiplication:
\begin{equation*}\begin{split}
xy&=(\alpha,\gamma_{12},\ldots, \gamma_{n-1,n})(\alpha',\gamma'_{12},\ldots,\gamma'_{n-1,n})\\
&=(\alpha+\alpha',
\gamma_{12}+\gamma'_{12}+f^i_{12}(\alpha,\alpha'),\ldots,\\
&\quad \gamma_{n-1,n}+\gamma'_{n-1,n}+f^i_{n-1,n}(\alpha,\alpha')),
\end{split}\end{equation*}
is assumed on $K_i$ and $f^i=(f^i_{12},\ldots,f^i_{n-1,n})$. Suppose $\eta_i:K_i\rightarrow
H_i$ be the group isomorphism whose existence established above. Now for each $1\leq i\leq n$ and $\alpha \in R$
let $h_i^\alpha\in H_i$ be the element of the equivalence class $(h_iZ(H))^\alpha$ such that
$$h_i^\alpha=\eta_i(\alpha,\underbrace{0,\ldots,0}_{(^n_2)-\textrm{times}}).$$
Firstly notice that $h_i^\alpha=1$ if and only if $\alpha=0$. Moreover it is clear that for each $1\leq i\leq n$
and $\alpha,\beta\in R$:
$$h_i^\alpha h_i^\beta=h_i^{\alpha+\beta}h_{12}^{f^i_{12}(\alpha,\beta)}\ldots
h_{n-1,n}^{f^i_{n-1,n}(\alpha,\beta)}.$$ Thus the relations (c) of Lemma~\ref{salam} hold between $h_i^\alpha$,
$1\leq i \leq n$ and $\alpha \in R$. We also note that,
\begin{equation*}\begin{split}[h_i^\alpha,h_j^\beta]&= f_H((h_iZ(H))^\alpha,(h_jZ(H))^\beta)\\
&=(f_H((h_iZ(H)),(h_jZ(H))))^{\alpha\beta}\\
&=[h_i,h_j]^{\alpha\beta}\\
&=h_{ij}^{\alpha\beta}\end{split}\end{equation*}
for $1\leq i<j\leq n$, and $\alpha,\beta\in R$, which proves that relations (a) hold.
Relations (b) hold since each $h_{ij}^\alpha$, $1\leq i <j \leq n$ and $\alpha\in R$, is central.
Relations (d) hold also in $H$ by the fact that each $H_{ij}$ is an $R$-module.

The set,
$$\mathcal{H}=\{h_i^\alpha,h_{kl}^\beta:1\leq i\leq n, 1\leq
k<l\leq n,\quad \alpha,\beta\in R \},$$ generates $H$ as a group by (5) of the definition of a basis. Let $F$ be the free group on $\mathcal{H}$.
Let $\mathcal{R}$ be the normal closure of the relations in the lemma~\ref{salam} in $F$, $g_i$ and $g_{ij}$
substituted by $h_i$ and $h_{ij}$ and the exponents come from the ring $R$ defined here and $\odot$ taken to be
concatenation. Let $\langle \mathcal{H}| \mathcal{R} \rangle$ be $F$ modulo the normal subgroup $\mathcal{R}$.
Consider the mapping:
$$\langle \mathcal{H}|\mathcal{R}\rangle \longrightarrow H, \quad h^\alpha_i \mapsto h_i^\alpha,
\quad h_{kl}^\alpha \mapsto h_{kl}^\alpha$$ for $\alpha \in R$, $1\leq i\leq n$ and $1\leq k<l \leq n$. The map
is a well defined homomorphism of groups since as proved above the relations (a)-(d) of~\ref{salam} hold also in
$H$. The map is also surjective since $\mathcal{H}$ generates $H$. Every word $W$ in $\langle
\mathcal{H}|\mathcal{R}\rangle$ is equivalent to a word of the form $h_n^{\alpha_n}\ldots
h_1^{\alpha_1}h_{12}^{\alpha_{12}}\ldots h_{n-1,n}^{\alpha_{n-1,n}}$. So $W$ maps to an element $h$ of $H$ of
this form. The uniqueness of this form for $h$ in $H$ is guarantied by (5) of~\ref{quasibasis}. So if $h$ is
trivial in $H$ then all the exponents in the above form are zero so the word $W$ is trivial in $\langle
\mathcal{H}|\mathcal{R}\rangle$. So $\langle \mathcal{H}|\mathcal{R}\rangle\cong H$. But by Lemma~\ref{salam},
$\langle \mathcal{H}|\mathcal{R}\rangle\cong N_{2,n}(R,f^1,\ldots f^n)$, hence
$$N^*_{2,n}(R,f^1,\ldots, f^n)\cong (H,\fra{b}).$$
\end{proof}
\begin{lem}\label{definability}Let $\fra{b}=\{h_1,\ldots, h_n\}$ be a basis for a group $H$.
Then the subgroups $Z(H)$, $H_i$, $1\leq i \leq n$, $H_{ij}$, $1\leq i<j \leq n,$ and $[H,H]$ are first order
definable in the enriched group $(H,\fra{b})$.
 Thus all the conditions (1)-(5) of the definition of basis can be expressed by first order formulas
 of an enriched signature of groups. Moreover the additional condition making a basis $\fra{b}$ a $P(f_H)$ basis is also expressible in the signature of groups.\end{lem}
\begin{proof}The center is defined by the formula
$$\phi_{Z(H)}(x):\forall y[x,y]=1.$$
Let $\bar{h}=(h_1, \ldots , h_n)$. Then For each $1\leq i\leq n$, $H_i$ is defined by:
$$\phi_{H_i}(\bar{h},x):[h_i,x]=1.$$
For each $1\leq i<j\leq n$, the subgroup $H_{ij}$ is generated by the set $\{[h_i,y]:y \in H_j\}$.
 So for every element $x$ of $H_{ij}$, for some fixed $1\leq i<j \leq n$, can be written as a product
 $$x=[h_i,y_1]\ldots [h_i,y_m],\quad y_1,\ldots, y_n \in H_j.$$
 Since $H$ is a 2-nilpotent group, by condition (4) of the definition of a basis, we can rewrite $x$ as $x=[h_i,y_1\ldots y_n]$.
 So the subgroup $H_{ij}$ is defined by the formula:
 $$\phi_{H_{ij}}(\bar{h},x):\exists y(x=[h_i,y]\wedge \phi_{H_j}(y)).$$
 By (4), $[H,H]$ sits inside the center. Therefore every element $x$ of $[H,H]$ has the form mentioned in
 (5)-(b). Conversely if some arbitrary element $x\in H$ has the form indicated in (5)-(b), then $x\in [H,H]$ since $H_{ij}\subseteq [H,H]$ for each $1\leq i<j\leq n$. Thus
 $[H,H]$ is defined in $(H,\mathcal{B})$ by the formula:
 $$\exists y_{12}\ldots y_{n-1,n}(\bigwedge_{1\leq i< j \leq n}\phi_{H_{ij}}(y_{ij})\wedge x=y_{12}\ldots y_{n-1,n}).$$
 It is now clear how to formulate conditions (1),(2),(3) and (5). Condition (4) is simply given by:
 $$\forall x,y,z([x,y].z=z.[x,y]).$$
Lemma~\ref{almost}shows clearly that there is a first order sentence $\phi$ of the signature of groups such that
$$H_i/Z(H)=(h_iZ(H))^{P(f_H)} \Leftrightarrow H\models \phi$$
for all $1\leq i \leq n$.
\end{proof}

\subsection{\protect Groups elementarily equivalent to a free 2-nilpotent group of arbitrary finite rank}

\begin{thm}\label{commut1} Let $R$ be a commutative associative ring with unit and $G=N_{2,n}(R,f_1,\ldots f_n)$.
 If a group $H$ is elementarily equivalent to $G$ then  $H\cong N_{2,n}(S,q_1, \ldots q_n)$
  for some ring $S$ such that $R\equiv S$. \end{thm}
 \begin{proof}The standard basis of $G$ is a $P(f_G)$-basis (see the (2) $\Rightarrow$ (1) of the proof of the characterization theorem). Since $G\equiv H$ the formulas that interpret $f_H$ in $H$ are the same as the ones which interpret $f_G$ in $G$. Moreover $f_H$ has the same width as $f_G$. So the formulas that interpret the action of $P(f_G)$ on $G/Z(G)$ in $G$ are the same as the ones that interpret the action of $P(f)H$ on $H/Z(H)$ in $H$. So $P(f_G)\equiv P(f_H)$ and moreover $H$ has to have a $P(f_H)$-basis. So by the characterization theorem
 $H\cong N_{2,n}(S,q^1, \ldots , q^n)$ for some ring
 $S=P(f_H)\equiv P(f_G)\cong R$
 and the relevant symmetric 2-cocycles $q^i$.\end{proof}

The main result of this paper is just a corollary of Theorem~\ref{commut1}
\begin{thm}\label{main}Let $G$ be a free 2-nilpotent group of rank $n$. If  $H$ is  a group elementarily equivalent to $G$ then $H$ has
the form $N_{2,n}(R,f_1,\ldots f_n)$ for some ring $R\equiv \mathbb{Z}$.\end{thm}
\begin{proof} By Proposition~\ref{freenil}, $G\cong N_{2,n}(\mathbb{Z})$. By Theorem ~\ref{commut1}, $H$ has the form
indicated in the statement of the theorem.\end{proof}

\section{Central extensions and elementary equivalence}\label{cse}
The aim of this section is to prove that for any two elementarily equivalent characteristic zero integral domains $R$ and $S$
$$ N_{2,n}(R,f^1,\ldots, f^n)\equiv N_{2,n}(S,g^1,\ldots, g^n)$$ for any symmetric 2-cocycles $f^1,\ldots f^n,g^1,\ldots g^n$ described before. Meanwhile we prove that any two central extensions of a torsion free abelian group $A$ by a torsion free abelian group $B$ are elementarily equivalent providing that the 2-cocycles corresponding to the extensions differ up to a 2-coboundary by a symmetric 2-cocycle.
\begin{lem}\label{isocent}
Let
$$0 \rightarrow A \rightarrow G \rightarrow B\rightarrow 0$$
be a central extension of the abelian groups $A$ and $B$. Let $(I,\mathcal{D})$ be an ultrafilter. Then $G^I/\mathcal{D}$ is isomorphic to a central extension of $A^I/\mathcal{D}$ by $B^I/\mathcal{D}$.\end{lem}
\begin{proof} Let $f$ be the 2-cocycle corresponding to the extension above. Let $[x]$ denote an element of $B^I/\mathcal{D}$ when $x\in B^I$. We define a 2-cocycle:
$$f^{\mathcal{D}}:B^I/\mathcal{D} \times  B^I/\mathcal{D} \rightarrow A^I/\mathcal{D}, \quad ([x],[y]) \mapsto [f^I(x,y)]$$
where$f^I:B^I\times B^I\rightarrow A^I$ is the 2-cocycle defined by $f^I(x,y)(i)=f(x(i),y(i))$, $i\in I$. To prove the well-definedness let $[x]=[z], [y]=[t] \in B^I/\mathcal{D}$. Then,
$A=\{i\in I: x(i)=z(i)\}\in \mathcal{D}$ and $B=\{i\in I:y(i)=t(i)\} \in \mathcal{D}$. If $C=A\cap B$ then $C\in \mathcal{D}$ and
$$C \subseteq \{i\in I:f(x(i),y(i))=f(z(i),t(i))\}=\{i\in I: f^I(x,y)(i)=f^I(z,t)(i)\},$$
which implies that $f^\mathcal{D}([x],[y])=f^{\mathcal{D}}([z],[t])$.

Let $H$ be the central extension of $A^I/\mathcal{D}$ by $B^I/\mathcal{D}$ induced by $f^{\mathcal{D}}$. we denote an element of $G^I/\mathcal{D}$ by $[(b,a)]$ and an element of $H$ by $([b],[a])$ where  $b\in B^I$ and $a\in A^I$.   Define:
$$\phi:G^I/\mathcal{D}\rightarrow H, \quad [(b,a)]\mapsto ([b],[a])$$
We prove that $\phi$ is an isomorphism of groups. To prove well-definedness let $[(b,a)]=[(d,c)]$ for $b,d \in B^I$ and $a,c\in A^I$. Let $$C=\{i\in I: (b(i),a(i))=(d(i),c(i))\}.$$ Then $ C\subseteq D,E$ where
$D=\{i\in I:b(i)=d(i)\}$ and $E=\{i\in I:a(i)=c(i)\}$. So $D,E \in \mathcal{D}$ since $C$ is. Hence $([b],[a])=([d],[c])$. Surjectivity of $\phi$ is clear. To prove injectivity assume that $([b],[a])=(0,0)$ for $b\in B^I$ and $a\in A^I$. Let $D'=\{i\in I:b(i)=0\}$ and $E'=\{i\in I:a(i)=0\}$. So $D', E' \in \mathcal{D}$. But
$$\{i\in I: (b(i),a(i))=(0,0)\}=D'\cap E'\in \mathcal{D},$$
which implies that $[(b,a)]=[(0,0)]$. It can be easily checked that $\phi$ is a homomorphism.
\end{proof}
Let $A$ and $B$ are two abelian groups and $f\in \textrm{Z}^2(B,A)$. By $(A,B,f)$ we denote the central extension of $A$ by $B$ corresponding to $f$.
\begin{lem} \label{elem-ext1}Let $A$ and $B$ be torsion free abelian groups and $G=(A,B,f)$ and $H=(A,B,g)$ be as above. If $f-g \in S^2(B,A)$ then $G\equiv H$.\end{lem}
\begin{proof}Let $(I,\mathcal{D})$ be an ultrafilter over which $A^I/\mathcal{D}$ is $\omega_1$-saturated.  By Lemma~\ref{isocent} there are $f^\mathcal{D}, g^{\mathcal{D}}\in H^2(B^I/\mathcal{D}, A^I/\mathcal{D})$ such that $G^I/\mathcal{D}\cong (A^I/\mathcal{D},B^I/\mathcal{D},f^\mathcal{D})$ and $H^I/\mathcal{D}\cong (A^I/\mathcal{D},B^I/\mathcal{D},g^{\mathcal{D}})$. By Theorem 1.11 of ~\cite{eklof},$A^I/\mathcal{D}$ is pure injective. Moreover $B^I/\mathcal{D}$ is torsion free since $B$ is. These facts imply that $$\textrm{Ext}(B^I/\mathcal{D},A^I/\mathcal{D}) =0.$$ Note also that $f^\mathcal{D}-g^{\mathcal{D}}\in S^2(B^I/\mathcal{D},A^I/\mathcal{D})$. Hence $f^\mathcal{D}$ and $g^{\mathcal{D}}$ are cohomologous. So
$$G^I/\mathcal{D}\cong (A^I/\mathcal{D},B^I/\mathcal{D},f^\mathcal{D}) \cong (A^I/\mathcal{D},B^I/\mathcal{D},g^{\mathcal{D}}) \cong H^I/\mathcal{D}.$$
Hence $G\equiv H$.\end{proof}
\begin{cor}\label{1234}Let $R$ be a characteristic zero integral domain and consider $N_{2,n}(R, f^1, \ldots , f^n)$ for some symmetric 2-cocycles $f^i$. Then
$$N_{2,n}(R,f^1, \ldots , f^n)\equiv N_{2,n}(R).$$\end{cor}
\begin{proof}Let
$$(\bar{\alpha},\bar{\beta}) =(\alpha_1, \ldots , \alpha_n, \beta_1, \ldots \beta_n)\in R^n\times R^n,$$
 and for each $1\leq i \leq n$ define
 $$q^i(\bar{\alpha}, \bar{\beta})=_{df}f^i(\alpha_i, \beta_i), \quad \forall (\bar{\alpha},\bar{\beta}) \in R^n\times R^n.$$
 We note that $N_{2,n}(R,f^1, \ldots , f^n)$ is defined by the 2-cocycle
 $$f=(f_{ij}):\oplus_{i=1}^nR^+ \times \oplus_{i=1}^nR^+ \rightarrow \oplus_{i=1}^{(^n_2)}R^+,$$
Where
$$f_{ij}(\bar{\al}, \bar{\be})=\al_i\be_j + \sum_{i=1}^nq^i(\bar{\al}, \bar{\be}), \quad \forall (\bar{\alpha},\bar{\beta}) \in R^n\times R^n.$$
Now to conclude we need to observe that $\sum_{i=1}^nq^i\in S^2(\oplus_{i=1}^n R^+,\oplus_{i=1}^{(^n_2)}R^+)$ and use Lemma~\ref{isocent}.\end{proof}

\begin{lem}\label{long}If $R\equiv S$ as rings then $N_{2,n}(R)\equiv N_{2,n}(S)$.\end{lem}
\begin{proof} We just need to make the simple observation that the group $N_{2,n}(R)$ is interpretable in the ring $R$ with the same formulas that interpret $N_{2,n}(S)$ in $S$. \end{proof}
The following statement is the obvious corollary of Lemma~\ref{1234} and Lemma~\ref{long}.
\begin{cor} \label{inverse}If $R$ and $S$ are elementarily equivalent characteristic zero integral domains then
$$N_{2,n}(R,f^1,\ldots, f^n)\equiv N_{2,n}(S,g^1,\ldots, g^n),$$
For any $f^i\in S^2(\oplus_{i=1}^n R^+,\oplus_{i=1}^{(^n_2)}R^+)$ and $g^i\in S^2(\oplus_{i=1}^nS^+,\oplus^{(^n_2)}_{i=1}S^+)$. \end{cor}

\section{A group elementarily equivalent to $N_{2,n}(\mathbb{Z})$ which is not $N_{2,n}$}\label{qnnotn}
In this section we prove the existence of a $QN_{2,n}$-group over a certain ring which is not a $N_{2,n}$-group over any commutative associative ring.
\begin{lem}\label{126a}
assume that $\varphi: G=N_{2,n}(R,f^1, \ldots , f^n)\rightarrow N_{2,n}(S)=H$ is an isomorphism of groups. If $\varphi_1:Ab(G)\rightarrow Ab(H)$ is the isomorphism induced by $\varphi$ and $\varphi_0:Z(G)\rightarrow Z(H)$ is the restriction of $\varphi$ to $Z(G)$ then there exists an isomorphism $\mu: R \rightarrow S$ such that
$$\varphi_1((xZ(G))^\alpha)=(\varphi(x)Z(H))^{\mu(\alpha)},\quad \textrm{for all  } x\in G \textrm{ and } \alpha\in R,$$
and
$$\varphi_0(x^\alpha)=(\varphi(x))^{\mu(\alpha)}, \quad \textrm{for all  } x\in Z(G) \textrm{ and } \alpha\in R.$$
\end{lem}
\begin{proof} Let $\theta: P(f_G) \rightarrow P(f_H)$ be the isomorphism obtained in Proposition~\ref{ringiso10}. By Proposition~\ref{P(f)-R} the map:
 $$\mu: R \rightarrow S, \quad \alpha_\phi \mapsto \alpha_{\theta(\phi)}$$ is an isomorphism of rings. Then

\begin{align*}
 \varphi_1((xZ(G))^{\alpha_\phi})&= \varphi_1\phi_1(xZ(G))\\
 &=\varphi_1\phi_1(\varphi_1)^{-1}\varphi_1(xZ(G))\\
 &=\theta(\phi_1)\varphi_1(xZ(G))\\
 &=(\varphi_1(xZ(G)))^{\mu(\alpha_\phi)}\\
 &=(\varphi(x)Z(H))^{\mu(\alpha_\phi)}
\end{align*}
 for all $x$ in $G$. A similar argument using $\phi_0$ instead of $\phi_1$ proves that
$$\varphi(x^\alpha)=\varphi(x)^{\mu(\alpha)}, \quad \textrm{for all  }x\in Z(G) \textrm{ and } \alpha \in R.$$
\end{proof}
\begin{thm}\label{importantcor}
Let $R$ be a binomial domain and $$\vphi:G=N_{2,n}(R,f^1, \ldots, f^n)\rightarrow N_{2,n}(S)=H$$ be an isomorphism of groups. Then for each $1\leq j \leq r$ we have that $f^j\in Z^2(R^+, \oplus_{i=1}^{\binom{n}{2}}R^+)$, i.e. each $f^j$ is a 2-coboundary.
\end{thm}
\begin{proof}
Let $\fra{b}=\{g_1, \ldots, g_n\}$ be a standard basis for $G$ and the $g_{ij}$ be defined as usual. Set $\vphi(g_i)=h_i$ and $\vphi(g_{ij})=h_{ij}$. Let $\vphi_1:Ab(G)\rightarrow Ab(H)$ be the group isomorphism induced by $\vphi$. By Lemma~\ref{126a} there exists an isomorphism $\mu:R \rightarrow S$ of rings so that $\vphi_1((xZ(G))^\alpha)=(\vphi_1(x\Gamma_2(G)))^{\mu(\alpha)}$, for all $\alpha$ in $R$ and $x$ in $G$. This implies that $\{h_1Z(H), \ldots , h_nZ(H)\}$ generates $H/\Gamma_2(H)$ freely as an $S$-module since $\{g_1Z(G), \ldots, g_nZ(G)\}$ generates $G/\Gamma_2(G)$ freely as an $R$-module. So $\fra{c}=\{h_1, \ldots , h_n\}$ generates $H$ as an $S$-group (see Remark~\ref{Rgr} and \cite{grun}, 4.1). So every element $h$ of $H$ has a unique representation
$$h= h_n^{\al_n}\cdots h_1^{\al_1}h_{12}^{\ga_{12}}\cdots h_{n-1,n}^{\ga_{n-1,n}},$$ where the $\al_i$ and $\ga_{ij}$ are elements of $S$. For simplicity let us denote by ${\bf h}_2^{\ga}$ the product $h_{12}^{\ga_{12}}\cdots h_{n-1,n}^{\ga_{n-1,n}}$. By Lemma~\ref{126a}
$$\phi(g_i^\alpha)=h_{i}^{\mu(\alpha)}\mathbf{h}_2^{g(\mu(\alpha))}, \quad \forall \alpha \in R,$$
where $g=(g_{ij}):S \rightarrow \prod_{i=1}^{\binom{n}{2}}S$ is a function determined by $\vphi$.
Choose two arbitrary elements $\beta$ and $\beta'$ in $S$. Then,

\begin{align*}
 h_i^{\beta+\beta'}&=h_i^\beta h_i^{\beta'}\\
 &=\vphi (g_i^{\mu^{-1}(\beta )}){\bf h}_2^{-g(\beta)}\vphi ((g_i)^{\mu^{-1}(\beta' )}){\bf h}_2^{-g(\beta')}\\
 &=\vphi (g_i^{\mu^{-1}(\beta )})\vphi (g_i^{\mu^{-1}(\beta' )}){\bf h}_2^{-g(\beta)-g(\beta')}\\
 &=\vphi (g_i^{\mu^{-1}(\beta +\beta')}{\bf h}_2^{f^i(\mu^{-1}(\beta ), \mu^{-1}(\beta'))}){\bf h}_2^{-g(\beta)-g(\beta')}\\
 &=h_i^{\beta+\beta'}{\bf h}_2^{\mu f^i (\mu^{-1}(\beta ), \mu^{-1}(\beta'))+g(\beta+\beta')-g(\beta)-g(\beta')}
 \end{align*}
The identity above clearly shows that the $\mu f^i (\mu^{-1}(-), \mu^{-1}(-))\in Z^2(S^+, \oplus_{i=1}^{n_c}S^+)$. Since $\mu$ is a ring isomorphism this implies that for each $1\leq j \leq r$, $f^i$ is a 2-coboundary as claimed.

\end{proof}
\begin{lem}[Belegradeck]There is a ring $R$, $R\equiv \mathbb{Z}$ such that $Ext(R^+,R^+)\neq 0$.\end{lem}
\begin{proof}See~\cite{beleg99}.\end{proof}
\begin{thm}Let $G\cong N_{2,n}(\mathbb{Z})$, i.e. $G$ be a free nilpotent group of rank $n$ and class 2. There is a $QN_{2,n}$ group $H$ such that $G\equiv H$ but $H$ is not an $N_{2,n}$ group over any commutative associative ring.\end{thm}
\begin{proof}Let $R$ be ring $R\equiv \mathbb{Z}$ such that $Ext(R^+,R^+)\neq0$. By Corollary~\ref{importantcor} there are 2-cocycles $f^i:R^+\times R^+\rightarrow \oplus^{(^n_2)}_{i=1} R^+$, $1\leq i \leq n$, such that
$$H=N_{2,n}(R,f^1,\ldots , f^n) \ncong N_{2,n}(R).$$

We note that $H\ncong N_{2,n}(S)$ for any associative commutative ring $S$ by Theorem~\ref{commut}. Moreover $H \equiv G$ by Lemma~\ref{inverse}.\end{proof}

\end{document}